\documentclass[pdflatex,sn-mathphys-num]{sn-jnl}

\usepackage{array}
\usepackage{graphicx}%
\usepackage{multirow}%
\usepackage{amsmath,amssymb,amsfonts}%
\usepackage{amsthm}%
\usepackage{mathrsfs}%
\usepackage[title]{appendix}%
\usepackage{xcolor}%
\usepackage{textcomp}%
\usepackage{manyfoot}%
\usepackage{booktabs}%
\usepackage{algorithm}%
\usepackage{algorithmicx}%
\usepackage{algpseudocode}%
\usepackage{listings}%

\begin{document}

\title[Sheets of Spectral Data of Stokes Waves in Weakly Nonlinear Models]{Sheets of Spectral Data of Stokes Waves in Weakly Nonlinear Models}


\author*[1]{\fnm{Benjamin} \sur{Akers}}\email{benjamin.akers@afit.edu}

\author[2]{\fnm{Ryan} \sur{Creedon}}\email{ryan\_creedon@brown.edu}
\equalcont{These authors contributed equally to this work.}

\affil*[1]{\orgdiv{Department of Mathematics and Statistics}, \orgname{Air Force Institute of Technology}, \orgaddress{\street{2950 Hobson Way}, \city{Wright Patterson AFB}, \postcode{45433}, \state{OH}, \country{USA}}}

\affil[2]{\orgdiv{Department of Applied Mathematics}, \orgname{Brown University}, \orgaddress{\street{75 Waterman St.}, \city{Providence}, \postcode{02912}, \state{RI}, \country{USA}}}



\abstract{We study the spectral stability of small-amplitude Stokes waves in a family of
weakly nonlinear, unidirectional models of the form
$u_t + L u + (u^2)_x = 0$. We introduce a perturbation method to expand
the spectral data in wave amplitude near flat-state eigenvalue collisions,
with the ratio of the colliding modes as a free parameter. This yields
sheets of spectral data whose slices at fixed amplitude give isolas of
instability. The same perturbation framework treats both high-frequency
and Benjamin--Feir instabilities, extends to discontinuous dispersion
relations (including the Akers--Milewski equation), and, for the first time, provides an analytic approximation of the Benjamin–Feir spectrum for this model and a
direct comparison of high-frequency and Benjamin--Feir growth rates across
the full family of models. Asymptotic predictions
are validated against numerical spectra computed by Floquet--Fourier--Hill and
quasi-Newton methods.}




\maketitle

\section{Introduction}

Periodic traveling water waves have been studied for almost two centuries, beginning with the work of George Stokes. Provided the amplitude $\varepsilon$ of these waves is sufficiently small, Stokes showed that these waves can be constructed as a power series in  $\varepsilon$. For this reason, we refer to these waves as \emph{Stokes waves}. A complete historical account of Stokes waves is given in~\cite{craik2005george}. 

The stability of Stokes waves has been
investigated with both multi-scale and direct spectral methods. Beginning in the
1960s, multi-scale expansions led to amplitude equations such as the nonlinear
Schr\"odinger, Dysthe, and Davey--Stewartson equations, in which the stability
of Stokes waves to modulational sideband perturbations was analyzed
in~\cite{benney1969wave,zakharov1968stability,phillips1967theoretical,benjamin1967disintegration,dysthe1979note,davey1974three,segur2005stabilizing}. For more general perturbations, one obtains
a spectral problem for the linearization of the governing water-wave equations
about a Stokes wave. The specta $\lambda$ are purely continuous and parameterized by
a Floquet exponent $p\in[0,1)$, encoding the quasi-periodicity of the
perturbations~\cite{kapitula2013spectral}. In the recent studies ~\cite{akers2015modulational,creedon2021higha,creedon2021highb,creedon2022high,creedon2023high}, asymptotic expansions of the unstable spectra and their corresponding
Floquet exponents were obtained for sufficiently small $\varepsilon$, revealing various instabilities of Stokes waves, including the classic \emph{modulational (or Benjamin--Feir) instability} as well as the newly discovered \emph{high-frequency instabilities}.  

In this work, we consider these asymptotic spectral expansions in the setting of a general
family of weakly nonlinear models of the form
\begin{equation}
  u_t + {L} u + (u^2)_x = 0,
  \label{mainModel}
\end{equation}
where ${L}$ is a Fourier multiplier with purely imaginary symbol
$\hat{{L}} = i \omega(k)$ and $\omega$ is the linear dispersion relation of the model.  Equations of the form~\eqref{mainModel}
are standard unidirectional models for surface water waves and related dispersive
systems. Concrete examples include the Whitham, capillary--Whitham, Kawahara,
and Akers--Milewski equations, whose dispersion relations we summarize in
Table~1. For each such model, small-amplitude Stokes waves exist and are
parametrically analytic in the amplitude $\varepsilon$~\cite{akers2021wiltonB,creedon2026existence}. We investigate the unstable spectra of these waves for $0 < \varepsilon \ll 1$, noting that, since the specific models under consideration are weakly nonlinear approximations of the full water wave equations, the instability mechanisms we find have direct implications for water-wave dynamics.

\begin{table}[tb]
\begin{tabular}{>{\centering\arraybackslash}m{3cm} >{\centering\arraybackslash}m{4cm}}
{\small \textbf{Model}} &{\small $\boldsymbol{\omega(k)}$ }\\ \hline\hline
\vspace*{0.1cm} Kawahara & \vspace*{0.1cm} \(a k^3 + b k^5\) \\
Whitham & \(k\sqrt{\frac{\tanh(kh)}{k}}\) \\
Capillary-Whitham & \(k\sqrt{\left(1 + \sigma k^2\right)\frac{\tanh(kh)}{k}}\) \\
Akers--Milewski \vspace*{0.2cm} & \(\textrm{sgn}(k)\left(1+\sigma |k|\right)^2\)  \vspace*{0.2cm} \\  \hline\hline 
\end{tabular} ~\vspace*{0.1cm}
\caption{ \it Example dispersion relations of weakly nonlinear models given by equation \eqref{mainModel}.  }
\label{tabOfModels}
\end{table}

Our first goal is to organize the spectral data near eigenvalue collisions of the
linearization about the flat state. At zero amplitude $(\varepsilon = 0)$, the
stability spectrum of~\eqref{mainModel} consists of imaginary eigenvalues
$\lambda_0(k,p)$, indexed by the integer wavenumber $k\in\mathbb{Z}$ and
Floquet exponent $p\in[0,1)$. We consider (\emph{i}) collisions of two nonzero
eigenvalues at a Floquet exponent $p_0\neq 0$, leading to
{high-frequency instabilities}, and (\emph{ii}) the collision of eigenvalues at
the origin with $p_0 = 0$, corresponding to the Benjamin--Feir instability.
For fixed $0 < \varepsilon \ll 1$, each collision of type (\emph{i}) gives rise to a curve of unstable spectra parameterized by an $\varepsilon$-dependent interval of Floquet exponents $p$. This curve, known as an \emph{isola of instability} in the water wave literature, is approximately an ellipse centered about the imaginary axis in the complex spectral plane. In contrast, for collision corresponding to the Benjamin-Feir instability, the curve is approximately a lemniscate centered at the origin. 

To obtain asymptotic expansions for these curves of unstable spectra, one expands the eigenvalues, eigenfunctions, and Floquet exponent as power series in $\varepsilon$. As the calculations unfold, one realizes that a coefficient of the Floquet exponent expansion can be chosen as a free parameter, which ultimately parameterizes the curve of unstable spectra for fixed $\varepsilon$. Varying both $\varepsilon$ and this parameter yields a 2D \emph{sheet of spectral data} in $(\textrm{Re} \lambda, \textrm{Im}\lambda,\varepsilon)$ space. 

Alternatively, one can choose the ratio $\beta$ of the amplitudes of the colliding flat-state eigenfunctions to parameterize the curve of unstable eigenvalues for fixed $\varepsilon$. We opt for this choice for two reasons. First, the solvability conditions that appear at each order in the perturbation scheme reduce to linear systems, simplifying the algebra of the expansion while
producing spectral data that agree with the previous parametrization after
relabeling. Second, varying $\beta$ naturally reveals important stability criterion like the classical Krein signature condition for high-frequency instabilities \cite{deconinck2017high}. In fact, a trivial consequence of our calculations shows that this Krein condition is both necessary and sufficient for the onset of high-frequency instabilities in this class of models. 

Building on earlier analyses of high-frequency and Benjamin--Feir instabilities
in specific models \cite{akers2015modulational, creedon2021higha, creedon2021highb, creedon2022high},
the present paper introduces three main novelties:
\begin{itemize}
  \item[1.] We develop a simple, direct perturbation method that produces both
  the high-frequency and Benjamin--Feir unstable spectra for the general
  model~\eqref{mainModel}. The method yields not only the leading-order
  behavior of the unstable eigenvalues, but also their associated eigenfunctions
  and the most unstable growth rates. For high-frequency instabilities arising
  from type (\emph{i}) collisions, we show that the structure of the resulting
  sheets of spectral data is organized by the resonant interaction between the
  carrier wave and the colliding modes: triad resonances ($|k_1-k_2|=1$) lead
  to elliptical isolas that are centred at the flat-state eigenvalue and have
  semi-major axes that scale as $\mathcal{O}(\varepsilon)$, while quartet
  resonances ($|k_1-k_2|=2$) produce elliptical isolas whose centres drift like
  $\mathcal{O}\left(\varepsilon^2\right)$ from the flat-state eigenvalue and
  have semi-major axes that scale as $\mathcal{O}\left(\varepsilon^2\right)$.
  \item[2.] Because the same perturbation framework applies to both the
  high-frequency and Benjamin--Feir instabilities across a family of weakly
  nonlinear models, we can directly compare the magnitude of their growth
  rates for different choices of the dispersion relation. In particular, for
  representative models such as the Whitham, capillary--Whitham, Kawahara,
  and Akers--Milewski equations, we quantify how the most unstable
  high-frequency eigenvalues compete with, or are dominated by, the
  Benjamin--Feir eigenvalues as the wave amplitude and wavenumber vary.
  \item[3.] We formulate the Benjamin--Feir perturbation method so that it
  remains valid even when the dispersion relation is not differentiable. In
  particular, we successfully obtain the Benjamin--Feir unstable eigenvalues of
  the Akers--Milewski equation, whose dispersion relation has a jump discontinuity at $k = 0$ and therefore falls outside the smooth setting
  of~\cite{creedon2023high} and related works. Our scheme only requires
  the pointwise values of $\omega(k)$ together with careful tracking of the
  singular terms that arise when $\omega$ has a jump. This allows us to
  obtain, to our knowledge, the first analytic approximation of the
  Benjamin--Feir spectrum for a model with discontinuous $\omega(k)$, and
  shows that the method we develop here is robust enough to treat equations
  with discontinuous dispersion, not just the classical smooth water-wave models.
\end{itemize}

Throughout, we compare the asymptotic predictions with numerical computations
of the spectrum using two methods that are independent of the asymptotics:
the Floquet--Fourier--Hill method~\cite{deconinck2006computing} and a quasi-Newton
continuation method for individual eigenpairs~\cite{claassen2018numerical}. For representative choices of the
dispersion relation---including the Whitham, capillary--Whitham, Kawahara,
Boussinesq--Whitham, and Akers--Milewski equations---we find asymptotically
correct agreement between the analytical and numerical spectra and highlight
examples where higher-order terms play a significant quantitative role.

The paper is organized as follows. In Section~2, we formulate the general model,
construct small-amplitude Stokes waves, and derive the associated spectral
problem for general quasi-periodic perturbations parameterized by a Floquet exponent. In Section~3, we introduce the
perturbation scheme for high-frequency instabilities near non-zero Floquet collisions and describe the resulting sheets
and isolas of spectrum. In Section~4, we adapt the same framework to the
Benjamin--Feir instability, derive the lemniscate approximating the figure-eight
spectrum, and extend the method to the Akers--Milewski equation with
discontinuous dispersion. We also compare directly the growth rates of the
high-frequency and Benjamin--Feir instabilities across our family of models.
Section~5 summarises our findings and outlines several directions for future
work, including rigorous existence results for unstable spectra based on the
bifurcation structure revealed here.

\section{Stokes Waves and Spectral Stability Problem}\label{sec:stokes-spectral}
\subsection{Expansion of Stokes Waves}
We seek periodic, traveling-wave solutions of~\eqref{mainModel} of the
form
\begin{equation}
  u(x,t) = u(\xi), \qquad \xi = x - c t, \label{travelingWave}
\end{equation}
where $c$ is the wave velocity. Without loss of generality, we restrict attention to $2\pi$–periodic
waves. Indeed, if one wishes to study a $2\pi/\kappa$–periodic traveling wave of
\eqref{mainModel} for any $\kappa > 0$, one can rescale the traveling wave coordinate by
\[
  \xi \mapsto \kappa^{-1} \xi,
\]
so that the wave becomes $2\pi$–periodic. Under this
rescaling, the dispersion relation and nonlinear term are modified in an elementary way, and $\kappa$ appears as an additional parameter in the model. Thus it
suffices to work with $2\pi$–periodic waves, with the understanding that
$\kappa$ can be reintroduced when desired.

Substituting \eqref{travelingWave} into
\eqref{mainModel} yields the steady traveling-wave equation
\begin{equation}
  -c u'(\xi) + L u(\xi) + \big(u^2(\xi)\big)' = 0,
  \label{eq:TW-equation}
\end{equation}
where primes denote differentiation with respect to $\xi$. For the remainder of this manuscript, we replace $\xi$ with $x$ for ease of notation. 

We are interested in small-amplitude solutions of \eqref{eq:TW-equation} that bifurcate from a single
Fourier mode. Let $\omega(k)$ denote the dispersion relation of~\eqref{mainModel}
and $c_p(k) = \omega(k)/k$ the corresponding phase velocity. For definiteness, we
take the carrier wavenumber to be $k=1$ so that the linear wave has phase
velocity
\[
  c_0 = c_p(1) = \omega(1).
\]
Following standard constructions (see, for example,
\cite{akers2021wiltonB}), we introduce a small amplitude parameter
$\varepsilon>0$ and seek solutions of~\eqref{eq:TW-equation} of the form
\begin{equation}
  u(x;\varepsilon) = \varepsilon \cos x
    + \varepsilon^2 u_2(x) + \varepsilon^3 u_3(x) + \cdots,
  \qquad
  c(\varepsilon) = c_0 + \varepsilon c_1 + \varepsilon^2 c_2 + \varepsilon^3 c_3 + \cdots,
  \label{eq:stokes-expansion-ansatz}
\end{equation}
with $2\pi$–periodic profiles $u_j$. Substituting~\eqref{eq:stokes-expansion-ansatz}
into~\eqref{eq:TW-equation} and equating powers of $\varepsilon$ yields a
hierarchy of linear problems for $u_j$, and solvability conditions of these linear problems determine the coefficients $c_j$.
To second order in $\varepsilon$, one finds
\begin{subequations}
\begin{align}
  u(x;\varepsilon)
    &= \varepsilon \cos x
      + \varepsilon^2 \tfrac{1}{2(c_0-c_p(2))}\cos(2x)
      + \mathcal{O}(\varepsilon^3),
 \\
  c(\varepsilon)
    &= c_p(1) + \varepsilon^2 \tfrac{1}{2(c_0-c_p(2))}
      + \mathcal{O}(\varepsilon^3),
\end{align}
\label{eq:stokes-expansion}
\end{subequations}
More generally, the Stokes expansion~\eqref{eq:stokes-expansion-ansatz} can be
continued to all orders in $\varepsilon$, and for each fixed dispersion relation
$\omega(k)$ the resulting power series converge for sufficiently small
$\varepsilon$; see~\cite{akers2021wiltonB,maspero2024full} for rigorous results in
this direction. In what follows, we regard $\varepsilon$ as the primary
bifurcation parameter for the family of Stokes waves and work with the expansion
\eqref{eq:stokes-expansion} as our model small-amplitude solution.

\subsection{Stability Spectrum of Stokes Waves}

To study the spectral stability of the Stokes wave $u(x;\varepsilon)$, we
consider perturbations of the form
\[
  u(x,t) = u(x;\varepsilon) + v(x,t),
\]
and linearize~\eqref{mainModel} about $u$. Writing $v(x,t)$ as a small
perturbation and discarding quadratic and higher order terms in $v$ yields the
linearized equation
\begin{equation}
  v_t - c(\varepsilon) v_x + {L} v + 2\big(u(x;\varepsilon)\,v\big)_x = 0.
  \label{eq:linearized-evolution}
\end{equation}
We seek normal modes of the form
\[
  v(x,t) = e^{\lambda t} v(x),
\]
where $\lambda\in\mathbb{C}$ is the spectral parameter and $v(x)$ satisfies a
Floquet (Bloch) boundary condition
\begin{equation}
  v(x + 2\pi) = e^{2\pi i p} v(x),
  \qquad p \in [0,1),
  \label{eq:floquet-condition}
\end{equation}
with Floquet exponent $p$ encoding the quasi-periodicity of the perturbation. Alternatively, we may write our normal mode as
\begin{align}
v(x,t) = e^{\lambda t}e^{i p x}\tilde{v}(x), \label{normalMode}
\end{align}
such that $\tilde{v}(x)$ is a $2\pi$-periodic function. 

Substituting \eqref{normalMode} into~\eqref{eq:linearized-evolution} and dropping the tilde on $v$ for ease of notation, we arrive at the
spectral stability problem
\begin{equation}
  \lambda v - c(\varepsilon) (\partial_x+ip) v + \mathcal{L} v(x)
    + 2(\partial_x+ip)\big(u(x;\varepsilon)\,v(x)\big) = 0, \quad \textrm{where} \quad \mathcal{L} := e^{-ipx}\circ L \circ e^{ipx}.
  \label{eq:spectral-problem}
\end{equation}
In Fourier space, the operator $\mathcal{L}$ can be better seen as a multiplier with symbol
\begin{equation}
    \hat{\mathcal{L}} = i\omega(k+p).
\end{equation}
For each fixed
$0 \leq \varepsilon \ll 1$ and $p\in[0,1)$, the spectrum of~\eqref{eq:spectral-problem}
consists of a discrete set of eigenvalues $\lambda(p;\varepsilon)$ provided $v(x) \in \textrm{L}^2(0,2\pi)$; see \cite{kapitula2013spectral} for more details.

At zero amplitude $(\varepsilon = 0)$, the Stokes wave reduces to the flat
state $u\equiv 0$ and the wave velocity $c(\varepsilon)$ reduces to $c_0$. In this
limit, the spectral problem~\eqref{eq:spectral-problem} simplifies to
\begin{equation}
  \lambda v(x) - c_0 v_x(x) + \mathcal{L} v(x) = 0.
  \label{eq:flat-spectral-problem}
\end{equation}
The eigenfunctions of~\eqref{eq:flat-spectral-problem} are
the Fourier modes
\[
  v_{k,p}(x) = e^{i(k+p)x}, \qquad k \in \mathbb{Z},
\]
and the corresponding eigenvalues are
\begin{equation}
  \lambda_0(k,p) = i\big(\omega(k+p) - c_0 (k+p)\big),
  \qquad k \in \mathbb{Z}, \quad p \in [0,1).
  \label{eq:flat-eigenvalues}
\end{equation}
The family $\{\lambda_0(k,p)\}_{k\in\mathbb{Z},\,p\in[0,1)}$ forms the
Floquet--Bloch spectrum of the flat state, consisting entirely of imaginary
eigenvalues.

As the amplitude $\varepsilon$ of the Stokes wave is increased from zero, the
spectrum of~\eqref{eq:spectral-problem} deviates from the flat-state spectrum
\eqref{eq:flat-eigenvalues}. Most importantly, bifurcations occur near repeated eigenvalues of
\eqref{eq:flat-eigenvalues}. For nonzero repeated eigenvalues, we have $(k_1,k_2,p_0)$ such that
\[
  \lambda_0(k_1,p_0) = \lambda_0(k_2,p_0) \neq 0,
\]
for some Floquet exponent $p_0 \neq 0$ and integers $k_1 \neq k_2$. These eigenvalues give rise to the high-frequency instabilities when $0 < \varepsilon \ll 1$. In addition, we have a repeated eigenvalue at $\lambda_0 = 0$ corresponding to the Benjamin--Feir instability. The perturbation schemes developed in the next two sections describe how
these bifurcations unfold and generate the high-frequency and
Benjamin--Feir instabilities summarized in the Introduction.

\section{Asymptotics of High Frequency Instabilities\label{HighFrequency}}
The unstable spectra approximated in this section are assumed to bifurcate from a flat-state ($\varepsilon=0$), repeated eigenvalue given by $\lambda_0 = \lambda(k_1,p_0)=\lambda(k_2,p_0) \neq 0$ for some $p_0\ne0$ and $k_1\ne k_2\in\mathbb{Z}$.  The corresponding eigenfunction can be written explicitly as
\begin{equation}\label{Eigenvector0}
 v_0=\exp(i(k_1+p_0)x)+\beta \exp(i(k_2+p_0)x)=\phi_1(x)+\beta \phi_2(x).
 \end{equation}
 Due to the repeated nature of $\lambda_0$, the eigenfunction $v_0$ can be an arbitrary superposition of the modes $\phi_1$ and $\phi_2$.  In what follows, we normalize $v_0$ so that the coefficient of $\phi_1$ is unity. The parameter $\beta$ is also assumed to be nonzero. Otherwise, one finds that the higher-order corrections to the eigenvalues are imaginary and do not lead to instability.
 
We continue the eigenpair $(\lambda_0,v_0)$ in $\varepsilon$ by assuming the following power series expansions:
\begin{align}
\lambda = \sum_{j = 0}^{\infty}\lambda_j \varepsilon^j \quad \textrm{and} \quad  v(\xi) &= \sum_{j = 0}^{\infty}v_j(\xi)\varepsilon^j, \label{lam-v-expansions}
\end{align}
where $(\lambda,v)$ is an eigenpair of \eqref{eq:spectral-problem}. We choose to normalize $v(\xi)$ so that $\left<v(\xi),\phi_1\right> = 1$, where $\left<\cdot,\cdot\right>$ is the standard inner-product on $\textrm{L}^2(0,2\pi):$
\begin{align}
\left<f,g \right> = \frac{1}{2\pi}\int_0^{2\pi} \overline{f}gdx \quad \textrm{for} \quad f,g \in \textrm{L}^2(0,2\pi).
\end{align}
Note that this choice of normalization implies $\left<v_0,\phi_1\right> = 1$, consistent with our normalization of $v_0$, as well as implies $\left<v_j,\phi_1\right> = 0$ for all $j \geq 2$.

In addition to the expansions \eqref{lam-v-expansions}, we allow the Floquet exponent to vary analytically in $\varepsilon$:
\begin{align}
p = \sum_{j=0}^{\infty} p_j\varepsilon^j,
\end{align}
as in ~\cite{akers2015modulational,creedon2021higha,creedon2021highb, creedon2022high}. By positing this expansion, we allow the interval of Floquet exponents parameterizing the unstable eigenvalues to vary with $\varepsilon$. Consequently, the operator $\mathcal{L}$ becomes $\varepsilon$-dependent and requires the following expansion about $\varepsilon = 0$:
\begin{align}
\hat{\mathcal{L}} &= i\omega\left(k+\sum_{j = 0}^{\infty}p_j\varepsilon^j\right) \nonumber \\
&= i\omega(k+p_0) + \varepsilon p_1c_g(k+p_0) + \varepsilon^2\left(p_2c_g(k+p_0) +\frac12p_1^2c_g'(k+p_0) \right) + \mathcal{O}\left(\varepsilon^3\right), \nonumber\\
&= \hat{\mathcal{L}}_0 + \varepsilon\hat{\mathcal{L}_1} + \varepsilon^2\hat{\mathcal{L}_2} + \mathcal{O}\left(\varepsilon^3\right). \label{Lj_defn}
\end{align}
Here, $c_g(k) = \omega'(k)$ is the group velocity and primes denote differentiation in $k$.

The $\mathcal{O}(\varepsilon)$ equation for the spectral data $(\lambda,v)$ is 
\begin{equation}\label{SecondSpectrum} \left(\lambda_0-c_0\partial_\xi+\mathcal{L}_0\right)v_1=-\lambda_1 v_0+c_0ip_1v_0-\mathcal{L}_1v_0-2(v_0u_1)_x.
\end{equation}
From the choice of eigenpair $(\lambda_0,v_0)$, the modes $\phi_1$ and $\phi_2$ are in the nullspace of the operator $(\lambda_0-c_0\partial_x + \mathcal{L}_0)$. In order for \eqref{SecondSpectrum} to admit a consistent solution for $v_1$, the Fredholm alternative enforces two solvability conditions, namely
\[ \langle \phi_j,-\lambda_1 v_0+c_0ip_1v_0-\mathcal{L}_1v_0-2(v_0u_1)_{x}\rangle =0, \quad \textrm{for} \quad j \in \{1,2\}.\]
These conditions reduce to
\begin{equation}\label{MatrixSolvability}
 \left(\begin{array}{cc} -1 & -i\tilde{c}_g(k_1+p_0)\\ -\beta & -i\tilde{c}_g(k_2+p_0)\beta\end{array}\right)\left(\begin{array}{c}\lambda_1\\ p_1\end{array}\right)=\left(\begin{array}{c}\langle \phi_1,2(v_0u_1)_{x}\rangle\\  \langle \phi_2,2(v_0u_1)_{x}\rangle\end{array}\right),
\end{equation}
 where we have introduced the group velocity in the traveling frame $\tilde{c}_g(k)=c_g(k)-c_0$ to simplify the appearance of these conditions. \\\\
\textbf{Remark 3.1.} To compute the corrections to the flat-state spectral data, one needs two bifurcation parameters to satisfy the two solvability conditions that appear at $\mathcal{O}\left(\varepsilon\right)$ and at subsequent orders. In this work, corrections to the Floquet exponent $p_j$ and the eigenvalue $\lambda_j$ are chosen to satisfy the solvability conditions, while
$\beta$ is a free choice. By varying $\beta$, one obtains sheets of spectral data. This is in contrast to the approach of ~\cite{akers2015modulational,creedon2021higha,creedon2021highb,creedon2022high}, where $\beta$ is expanded as a power series in $\varepsilon$ and its corrections, alongisde the eigenvalue corrections, are the solvability parameters.  The corrections to the Floquet exponent in that framework are free. \\\\
\textbf{Remark 3.2.} Regardless of the frequency difference in the eigenfunctions $k_2 - k_1$, one always solves a linear system for $p_1$ and $\lambda_1$.  This contrasts \cite{akers2015modulational,creedon2021higha,creedon2022high,creedon2021highb}, where the solvability conditions are nonlinear.  Moreover the matrix which must be inverted for these (and the later) solvability conditions has a clear nondegeneracy condition given by the determinant
\[ \text{det}\left(\begin{array}{cc} -1 & -i\tilde{c}_g(k_1+p_0)\\ -\beta & -i\tilde{c}_g(k_2+p_0)\beta\end{array}\right)= i\beta( \tilde{c}_g(k_2+p_0)-\tilde{c}_g(k_1+p_0)). \]
If $\tilde{c}_g(k_2+p_0)-\tilde{c}_g(k_1+p_0)\ne 0$ then the solution formally exists for all $\beta\ne 0$. \\\\
\textbf{Remark 3.3.} Equation \eqref{MatrixSolvability} is forced by nonlinear interactions
\begin{equation}\label{InnerProd}
 \langle \phi_j,2(v_0u_1)_{x}\rangle=\frac{1}{2\pi}\int_0^{2\pi} e^{-i(k_j+p_0)x}\Big((e^{ix}+e^{-ix})(e^{i (k_1+p_0)x}+\beta e^{i(k_2+p_0)x}   )\Big)_{x} dx.
\end{equation}
For the above integrals to be nonzero, one needs a relationship between the wavenumbers of the eigenfunctions, $k_1$ and $k_2$, and that of the carrier wave, $k=1$.  The relationship is a triad condition, $k_1=k_2\pm 1$, without which there is no instability at this order.  Later orders will have a similar requirement, where the solvability conditions give rise to homogeneous linear equations unless the wavenumbers take part in a resonance interaction.  

Before proceeding with our analysis of \eqref{MatrixSolvability}, we note that the repeated eigenvalue $\lambda_0$ implies a resonance between the wavenumbers of the eigenmodes $k_j$ and that of the carrier wave $k=1$,  namely
\[ \lambda(k_1,p_0)= \lambda(k_2,p_0) \implies ic_0(k_1+p_0)-i\omega(k_1+p_0)=ic_0(k_2+p_0)-i\omega(k_2+p_0).\]
Using $c_0=c_p(1)=\omega(1)$ and defining $m=k_1-k_2$, we can rewrite this resonance condition as  
\begin{align} m-k_1+k_2=0\quad \textrm{and} \quad  m\omega(1)-\omega(k_1+p)+\omega(k_2+p)=0,\label{res-coll-cond}
\end{align}
which represents a resonance between $|m|+2$ waves. For $|m|>1$, the resonance is degenerate because the wavenumber $k=1$ is repeated $|m|$ times.  In the language of resonant interaction theory \cite{craik1988wave,hendersonRIT}, the $|m| = 1$ resonance represents a triad interaction, the $|m| = 2$ resonance represents a quartet interaction, and the $|m| \geq 3$ resonances represent high-order interactions. As will be seen, each value of $m$ gives rise to a different high-frequency instability as $\varepsilon$ increases from zero. Therefore, it is natural for us to label high-frequency instabilities according to the value of $m$ in \eqref{res-coll-cond} as triad ($|m| = 1$), quartet ($|m| = 2$), and higher-order ($|m| \geq 3$) high-frequency instabilities.  

\subsection{High Frequency Instabilities: Triads}

In this subsection, we assume that the wavenumbers $k_j$ of the flat-state eigenfunction take part in a triad resonance so that $|k_1-k_2|=1$.  For notational convenience we will define $k_2$ to be the larger of the two wavenumbers so that $k_2=k_1+1$.  No generality is lost here: since $\beta$ is already restricted to not be zero, swapping the roles of $k_1$ and $k_2$ can be accomplished by dividing $v$ by $\beta$.  Having chosen the triad resonance for $k_1$ and $k_2$, the right hand side of equation \eqref{MatrixSolvability} can be evaluated explicitly, yielding
\begin{equation}\label{MatrixSolveTriad}
 \left(\begin{array}{cc} -1 & -i\tilde{c}_g(k_1+p_0))\\ -\beta & -i\tilde{c}_g(k_2+p_0))\beta\end{array}\right)\left(\begin{array}{c}\lambda_1\\ p_1\end{array}\right)=\left(\begin{array}{c}  \beta i(k_1+p_0) \\   i(k_2+p_0) \end{array}\right),
\end{equation}
whose solution is
\[ \lambda_1=\frac{i\left\{\tilde{c}_g(k_1+p_0))(k_2+p_0)-\beta^2(k_1+p_0)\tilde{c}_g(k_2+p_0)\right\}}{\beta(\tilde{c}_g(k_2+p_0)-\tilde{c}_g(k_1+p_0))}, \]
\[ p_1=\frac{\beta^2(k_1+p_0)-(k_2+p_0)}{\beta( \tilde{c}_g(k_2+p_0)-\tilde{c}_g(k_1+p_0))}= \frac{\beta(k_1+p_0)- (k_2+p_0)\beta^{-1}}{( \tilde{c}_g(k_2+p_0)-\tilde{c}_g(k_1+p_0))}. \]
The eigenvalue correction, $\lambda_1$, is clearly imaginary for real $\beta$.   However, since the eigenfunction $v$ of this spectral problem can be complex-valued, $\beta$ can be complex. In contrast, the Floquet exponent $p$ must be real, so $p_1\in \mathbb{R}$. Writing $\beta=\rho \exp(i\theta,)$
\[ p_1= \frac{\rho (\cos(\theta)+i\sin(\theta))(k_1+p_0)- (k_2+p_0)\rho^{-1}(\cos(\theta)-i\sin(\theta))}{( c_g(k_2+p_0)-c_g(k_1+p_0))}. \]
To have real $p_1$, one requires
\[ (\rho(k_1+p_0)+(k_2+p_0)\rho^{-1})\sin(\theta) =0.\]
Thus, either $\sin(\theta)=0$ (and $\beta$ is real) or $\rho^2=\frac{-(k_2+p_0)}{(k_1+p_0)}$ has real solutions.  Since $\beta \in \mathbb{R}$ gives a purely imaginary (and, hence, stable) eigenvalue correction and since $\rho$ is real by definition, one needs
\[ \frac{k_2+p_0}{k_1+p_0}= \frac{k_1+1+p_0}{k_1+p_0}<0,\]
to achieve instability. The above condition is, in fact, a representation of the Krein signature computed in \cite{deconinck2017high,creedon2021higha}, which is a necessary condition for the formation of high-frequency instabilities. As we have just shown, it is also a sufficient condition. 

\begin{figure}[tp]
\centerline{\includegraphics[width=0.5\textwidth]{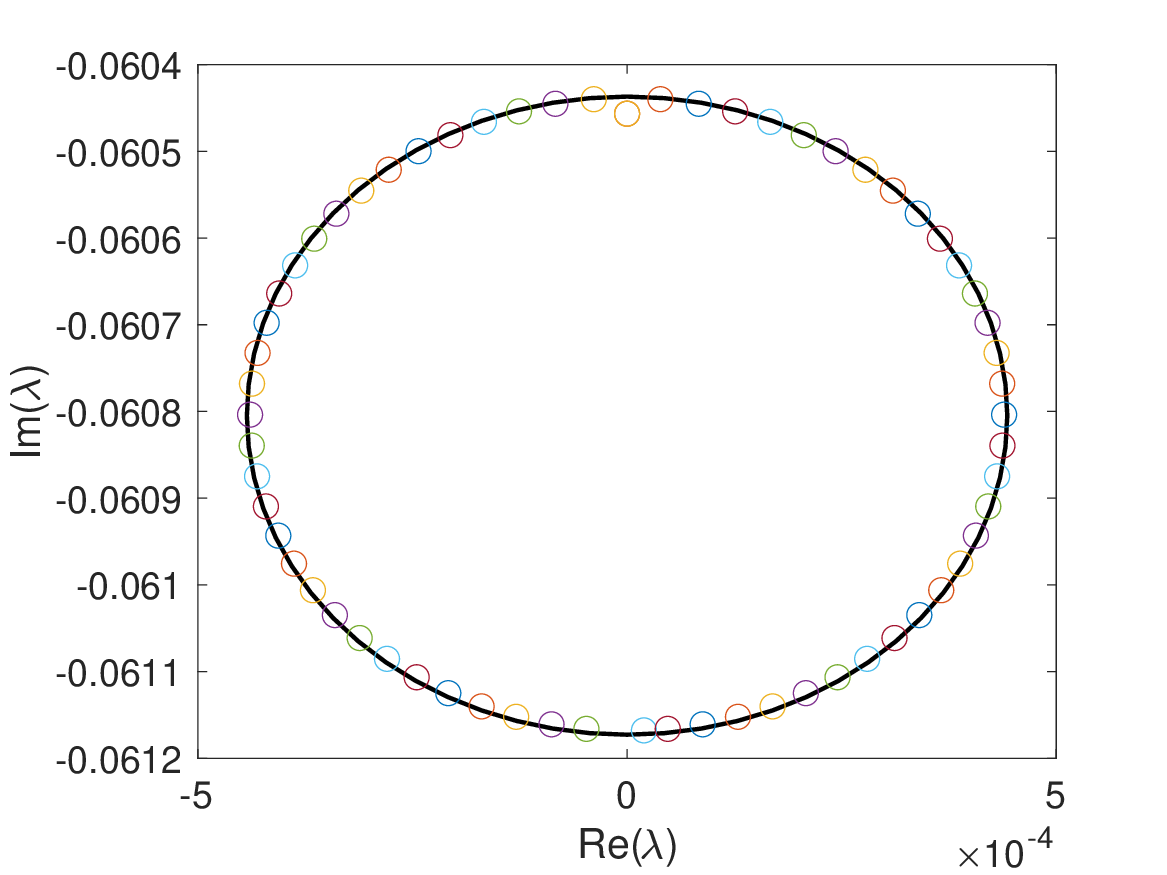}\includegraphics[width=0.5\textwidth]{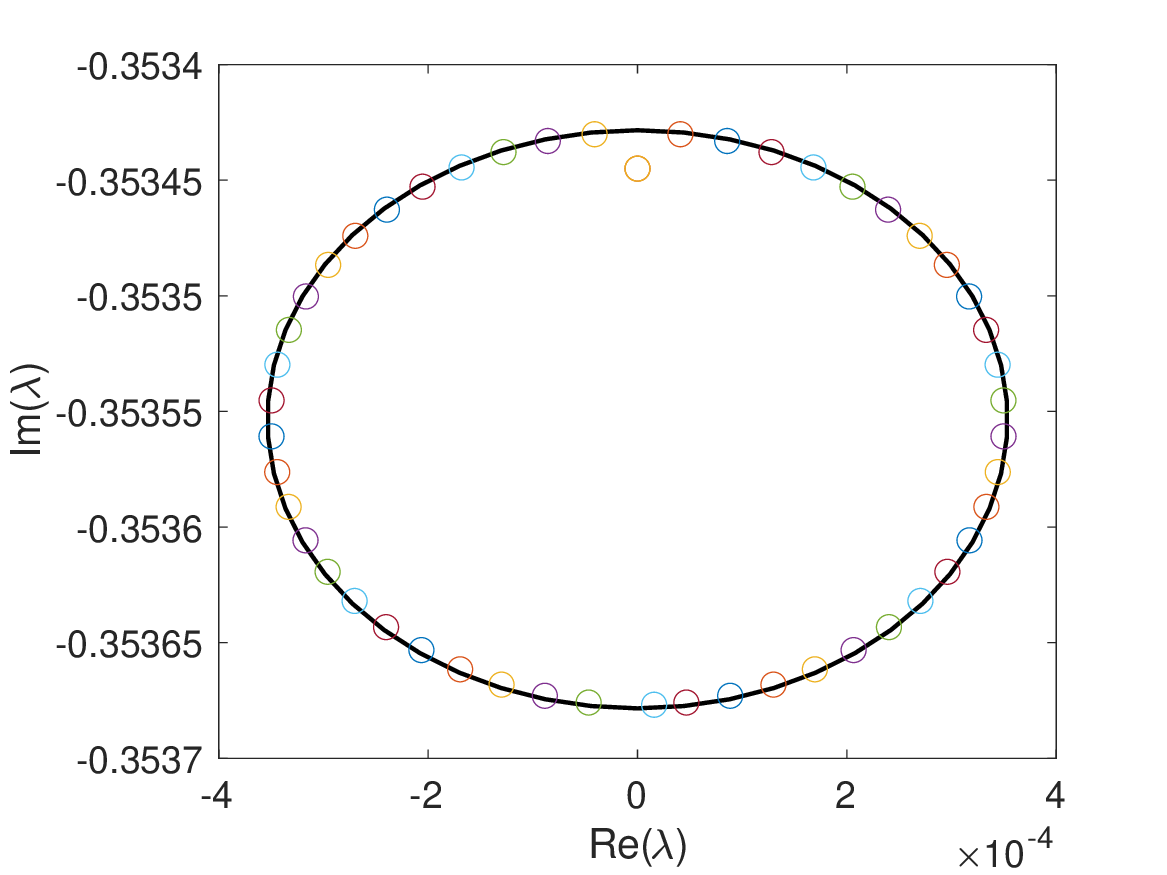}}
\caption{\it \textbf{Left:} An isola of a triad instability in the deep water ($h\rightarrow\infty$) gravity-capillary Whitham equation near $(\lambda_0,p_0)\approx(-0.0608i,0.2681)$  with $\sigma=2.5$ at $\varepsilon=10^{-3}$. The circles are computed with a quasi-Newton iteration; the curves are the asymptotic prediction of \eqref{TriadLam1}. \textbf{Right:} An isola of a triad instability in the Akers--Milewski equation near $(\lambda_0,p_0)=(0.3536i,0.1464)$  with $\sigma=2$ at $\varepsilon=10^{-3}$. The circles are computed with a quasi-Newton iteration; the curves are the asymptotic prediction of \eqref{TriadLam1}. \label{TriadIsolaFig}}
\end{figure}

\begin{figure}[tp]
\centerline{\includegraphics[width=0.5\textwidth]{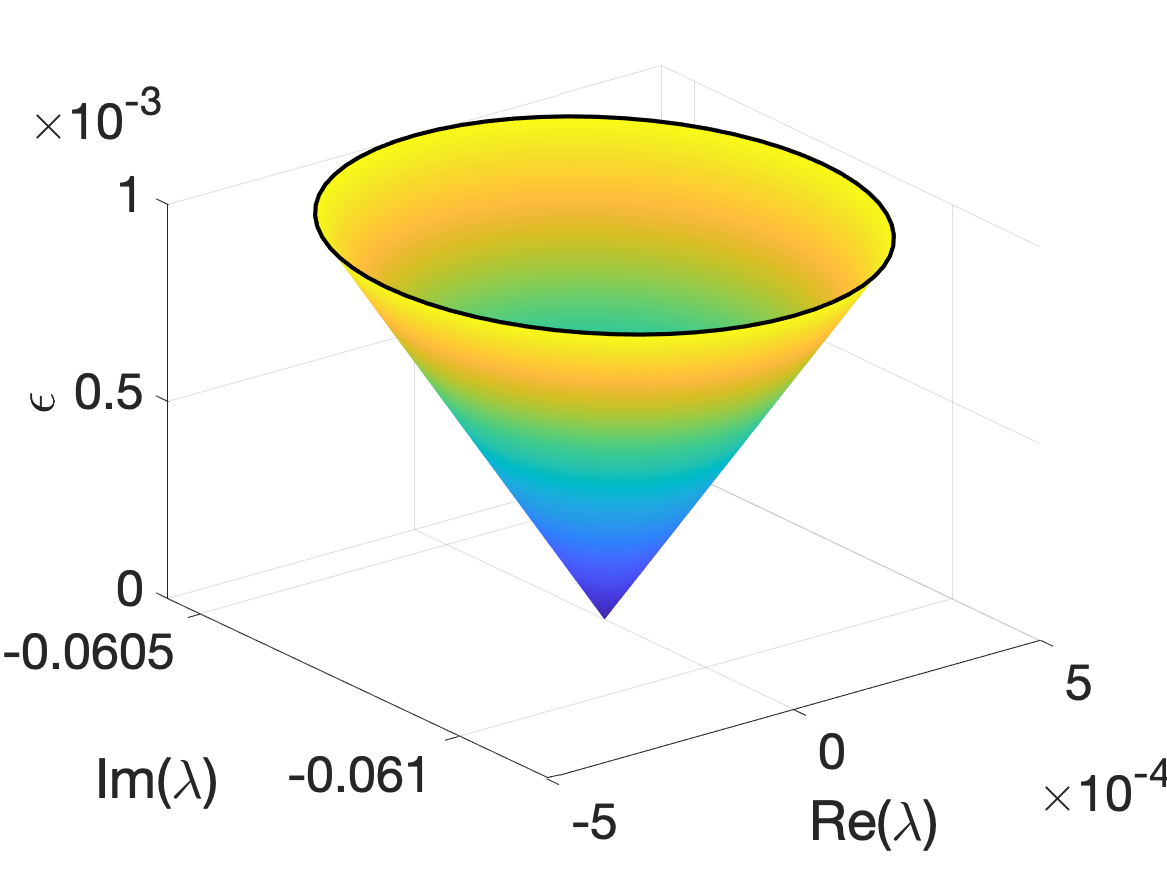}}
\caption{\it The surface of spectral data predicted by \eqref{TriadLam1} in the deep water ($h\rightarrow\infty$) gravity-capillary Whitham equation is visualized.  The spectrum bifurcated from $(\lambda_0,p_0)\approx(-0.0608i,0.2681)$ with $\sigma=2.5$.  This is the same configuration which is compared to numerical predictions in Figure \ref{TriadIsolaFig}.\label{TriadSurface}}
\end{figure}

Without loss of generality, we choose $k_1=-1$, $k_2=0$, and $\rho=\sqrt{\frac{p_0}{(1-p_0)}}$, giving us
\[ p_1= \left(\frac{\rho (k_1+p_0)- (k_2+p_0)\rho^{-1}}{( \tilde{c}_g(k_2+p_0)-\tilde{c}_g(k_1+p_0))}\right)\cos(\theta), \]
for any $\theta\in (0,2\pi)$ and 
\begin{subequations}\label{TriadLam1}
\begin{eqnarray}
 \text{Im}(\lambda_1)&=&\frac{\left(-\rho(k_1+p_0)\tilde{c}_g(k_2+p_0)+\rho^{-1}\tilde{c}_g(k_1+p_0))(k_2+p_0)\right)}{(\tilde{c}_g(k_2+p_0)-\tilde{c}_g(k_1+p_0))}\cos(\theta), \\
\text{Re}(\lambda_1)&=&-\frac{\left(\rho(k_1+p_0)\tilde{c}_g(k_2+p_0)+\rho^{-1}\tilde{c}_g(k_1+p_0))(k_2+p_0)\right)}{(\tilde{c}_g(k_2+p_0)-\tilde{c}_g(k_1+p_0))}\sin(\theta).
 \end{eqnarray}
 \end{subequations}
Thus, $\lambda_1$ predicts an elliptical isola centered at $\lambda_0$ for any sufficiently small $\varepsilon$. This ellipse is parameterized by $\theta$. Eigenvalues on this isola with positive real part are responsible for generating triad high-frequency instabilities. Example isolas for the deep water gravity-capillary Whitham equation and the Akers--Milewski equation are shown in Figure \ref{TriadIsolaFig}.  When both $\varepsilon$ and $\theta$ are varied, \eqref{TriadLam1} predicts a sheet of spectral data arranged on a conical surface; see, for instance, Figure \ref{TriadSurface}.

\subsubsection{High Frequency Instabilities: Quartets}
If the wavenumber difference is larger than one, $|k_1-k_2|>1$, then the integrals in \eqref{InnerProd} vanish, and the corrections $p_1$ and $\lambda_1$ are zero. Thus, there is no instability at $\mathcal{O}\left(\varepsilon\right)$.  Equation \eqref{SecondSpectrum} can be readily solved for $v_1$, giving us
\[ v_1=\alpha_{1}^{-} e^{i(k_1+p_0-1)x}+\alpha_{1}^{+} e^{i(k_1+p_0+1)x}+\beta\alpha_{2}^{-} e^{i(k_2+p_0-1)x}+\beta\alpha_{2}^{+}e^{i(k_2+p_0+1)x},\]
where $\alpha_j^{\pm}=\frac{-i(k_j+p_0\pm 1)}{\lambda_0(k_j+p_0)-\lambda_0(k_j+p_0\pm 1)}$ are real coefficients.\\

At $\mathcal{O}(\varepsilon^2)$, the spectral problem \eqref{eq:spectral-problem} becomes 
\begin{equation}\label{ThirdSpectrum} (\lambda_0-c_0\partial_{x}+\mathcal{L}_0)v_2=-\lambda_2 v_0+c_2v_{0,x}+c_0ip_2v_0-\mathcal{L}_1v_1 -\mathcal{L}_2v_0-2(v_1u_1+v_0u_2)_x,
\end{equation}
where $\mathcal{L}_1$ and $\mathcal{L}_2$ are Fourier multipliers with symbols defined in \eqref{Lj_defn}. Following the same calculations for the solvability conditions as in the $\mathcal{O}(\varepsilon)$ problem, one arrives at
\begin{equation}\label{MatrixSolvabilityQuartet}
 \left(\begin{array}{cc} -1 & -i\tilde{c}_g(k_1+p_0)\\ -\beta & -i\tilde{c}_g(k_2+p_0)\beta\end{array}\right)\left(\begin{array}{c}\lambda_2\\ p_2\end{array}\right)=\left(\begin{array}{c}\langle \phi_1,2(v_0u_2+v_1u_1)_{x}\rangle-c_2i(k_1+p_0)\\  \langle \phi_2,2(v_0u_2+v_1u_1)_{x}\rangle-c_2i(k_2+p_0)\beta\end{array}\right)
\end{equation}
The coefficient matrix above is identical to the triad setting and invertible whenever $\beta\ne0$ and $\tilde{c}_g(k_1+p_0)\ne \tilde{c}_g(k_1+p_0)$.  The right-hand side represents new possible resonances between $k_1$ and $k_2$. 
\subsubsection{Quartet Resonances}
\begin{figure}[tp]
\centerline{\includegraphics[width=0.5\textwidth]{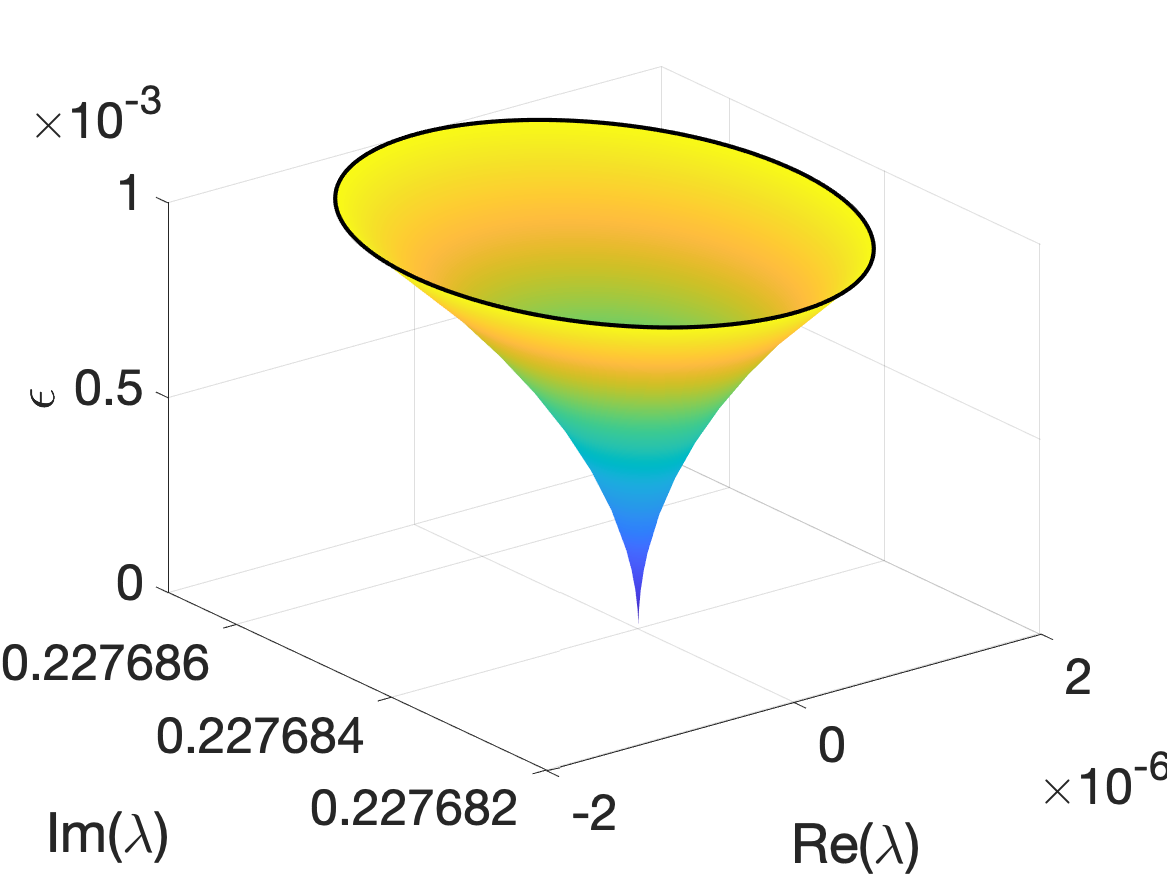}\includegraphics[width=0.5\textwidth]{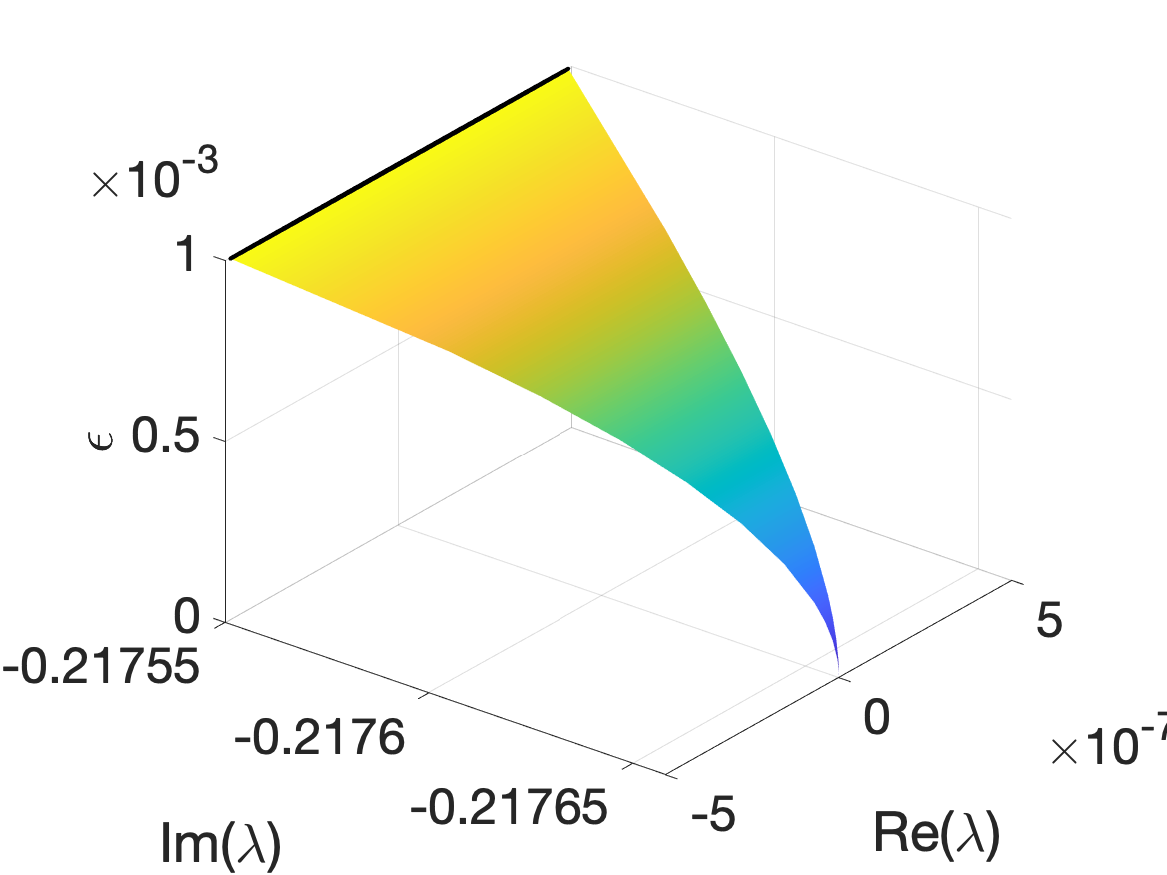}}
\caption{\it \textbf{Left:}The surface of spectral data predicted by \eqref{QuartetLam2} in the Kawahara equation is visualized.  The spectrum bifurcates from $(\lambda_0,p_0)\approx(0.2277i,0.3675)$ with $(a,b)=(1,-0.25)$.  This is the same configuration which is compared to numerical predictions in Figure \ref{KawaharaQuartet}.  \textbf{Right} The surface of spectral data predicted by \eqref{QuartetLam2} in the deep water gravity-capillary Whitham equation is visualized.  The spectrum bifurcates from $(\lambda_0,p_0)\approx(0.2177i,0.1363)$  with $\sigma=0.25$.  This is the same configuration which is compared to numerical predictions in Figure \ref{WhithamQuartet}.\label{QuartetSurface}}
\end{figure}

\begin{figure}[tp]
\centerline{\includegraphics[width=0.35\textwidth]{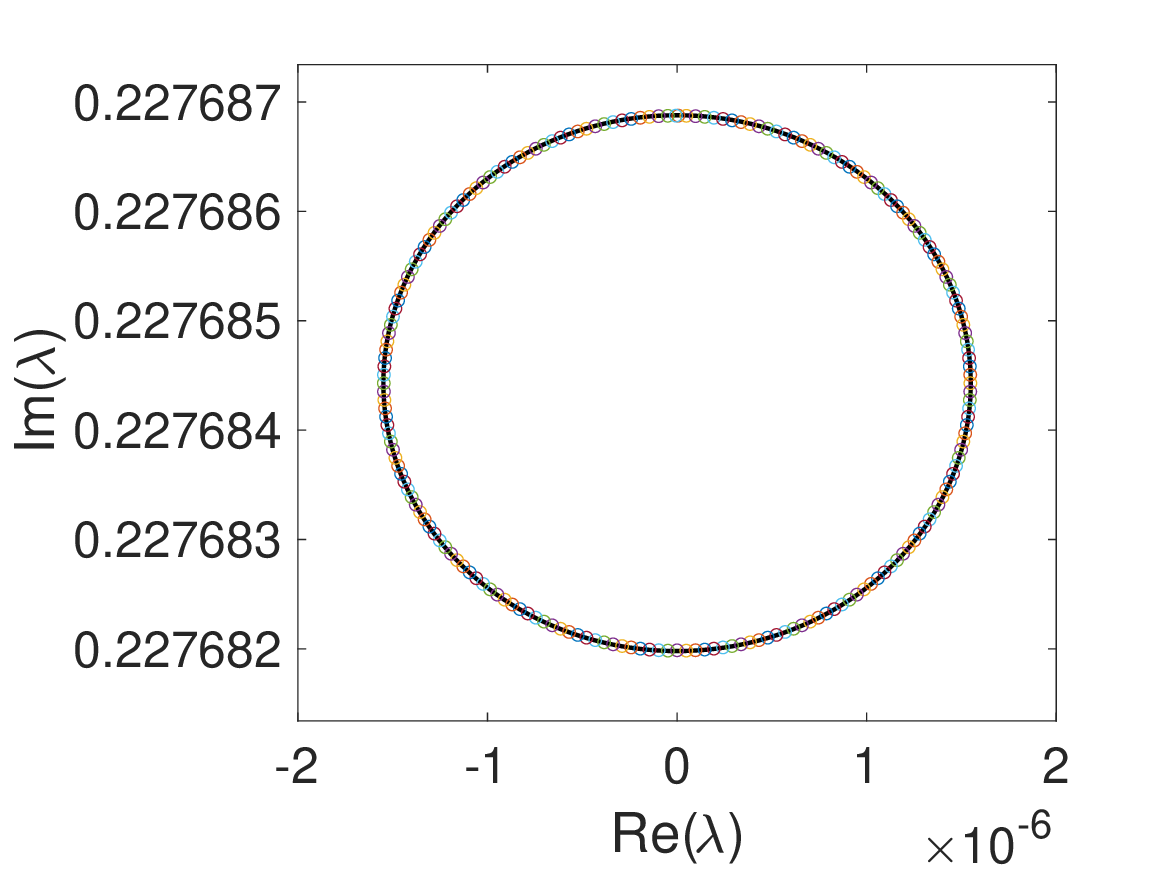}\includegraphics[width=0.35\textwidth]{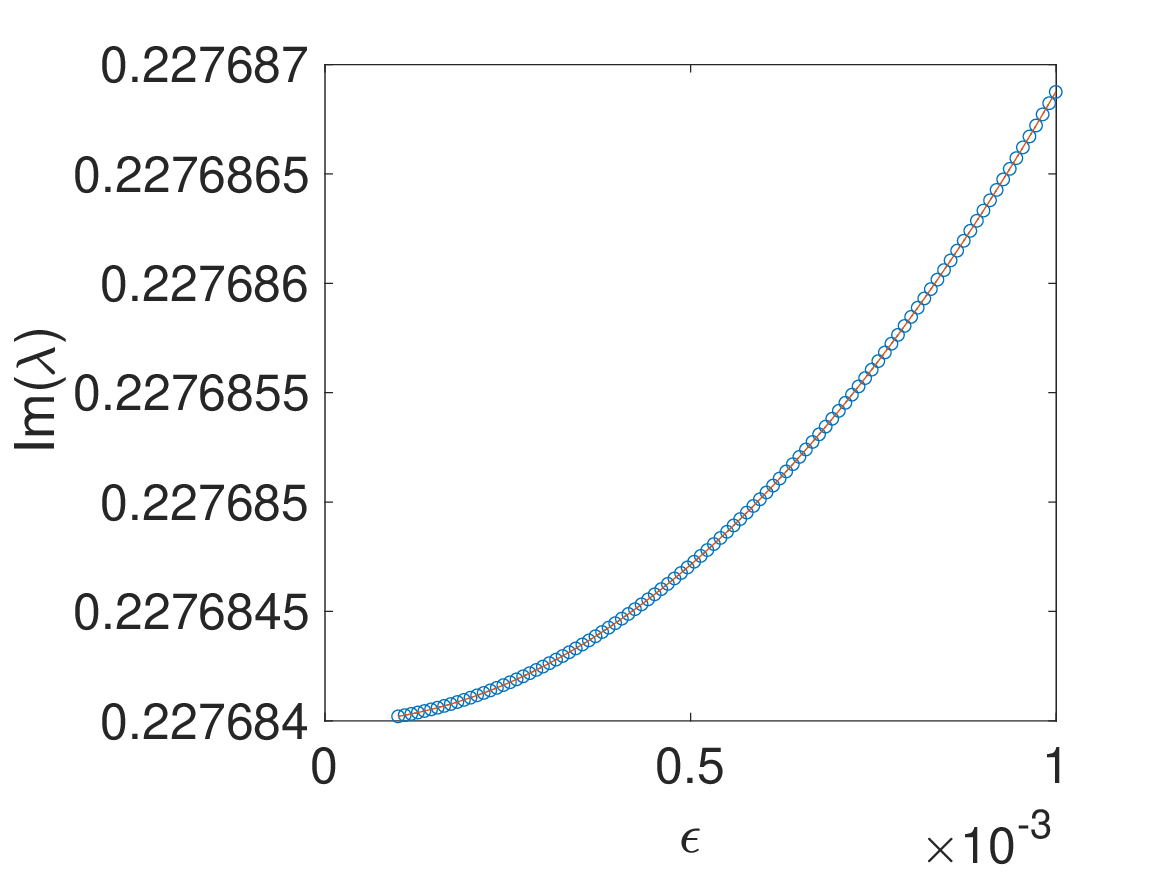}\includegraphics[width=0.35\textwidth]{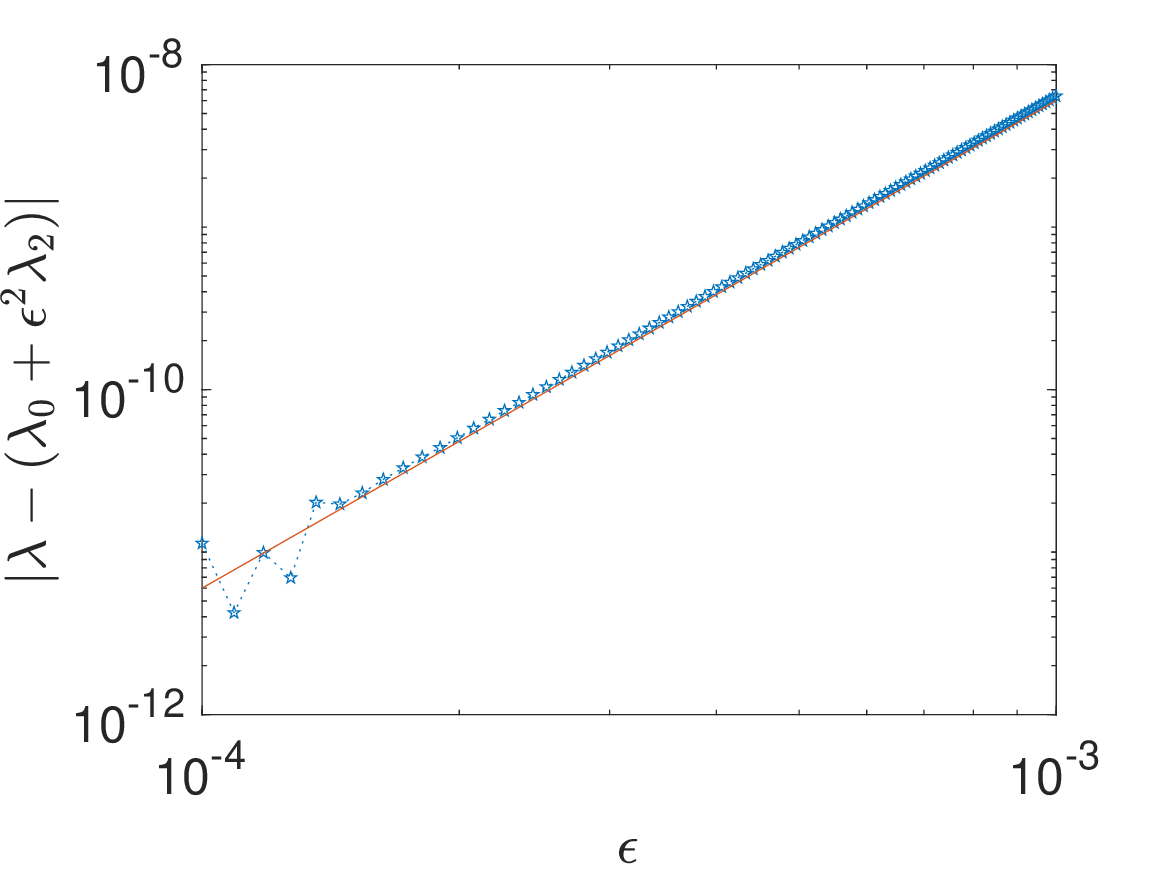}}
\caption{\it \textbf{Left:} An isola of a quartet instability in the Kawahara equation near $(\lambda_0,p_0)=(0.2277i,0.3675)$  with $(a,b)=(1,-0.25)$ at $\varepsilon=10^{-3}$. The circles are computed with a quasi-Newton iteration; the curves are the second order asymptotic prediction of \eqref{QuartetLam2}. \textbf{Center:} The imaginary part of the most unstable eigenvalue of \eqref{QuartetLam2} (solid line) against numerical computations from a quasi-Newton iteration (circles).  \textbf{Right:} The error between the numerically computed and asymptotics \eqref{QuartetLam2} (starred line); the solid line marks $y=6\varepsilon^3$.  \label{KawaharaQuartet}}
\end{figure}

\begin{figure}[htp]
\centerline{\includegraphics[width=0.35\textwidth]{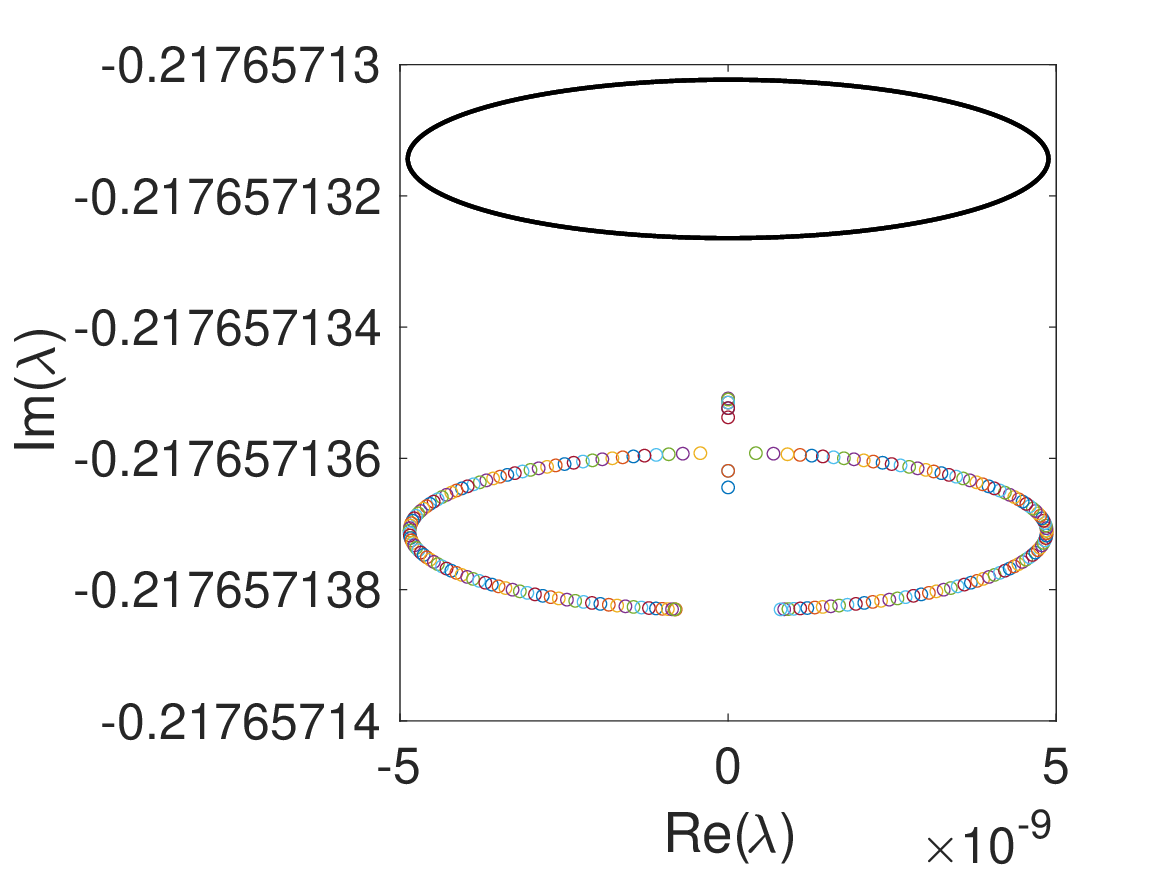}\includegraphics[width=0.35\textwidth]{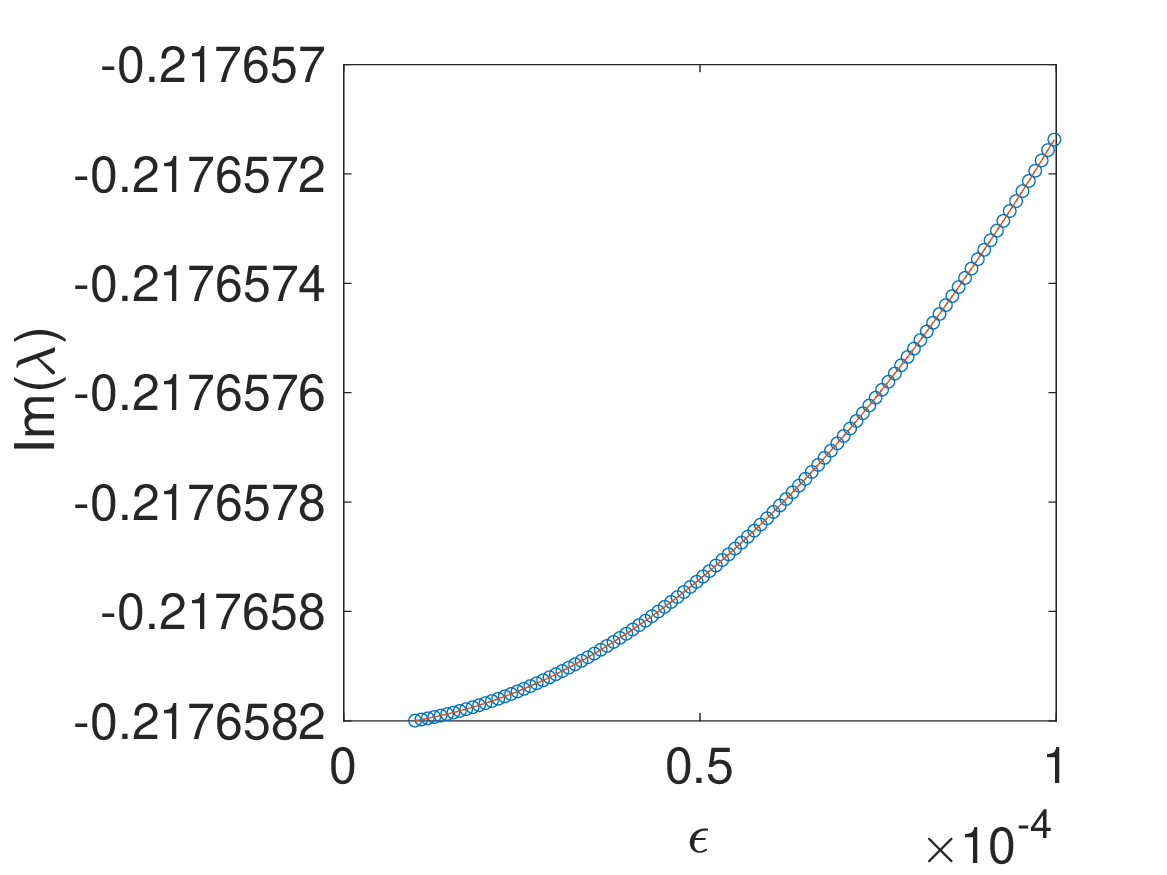}\includegraphics[width=0.35\textwidth]{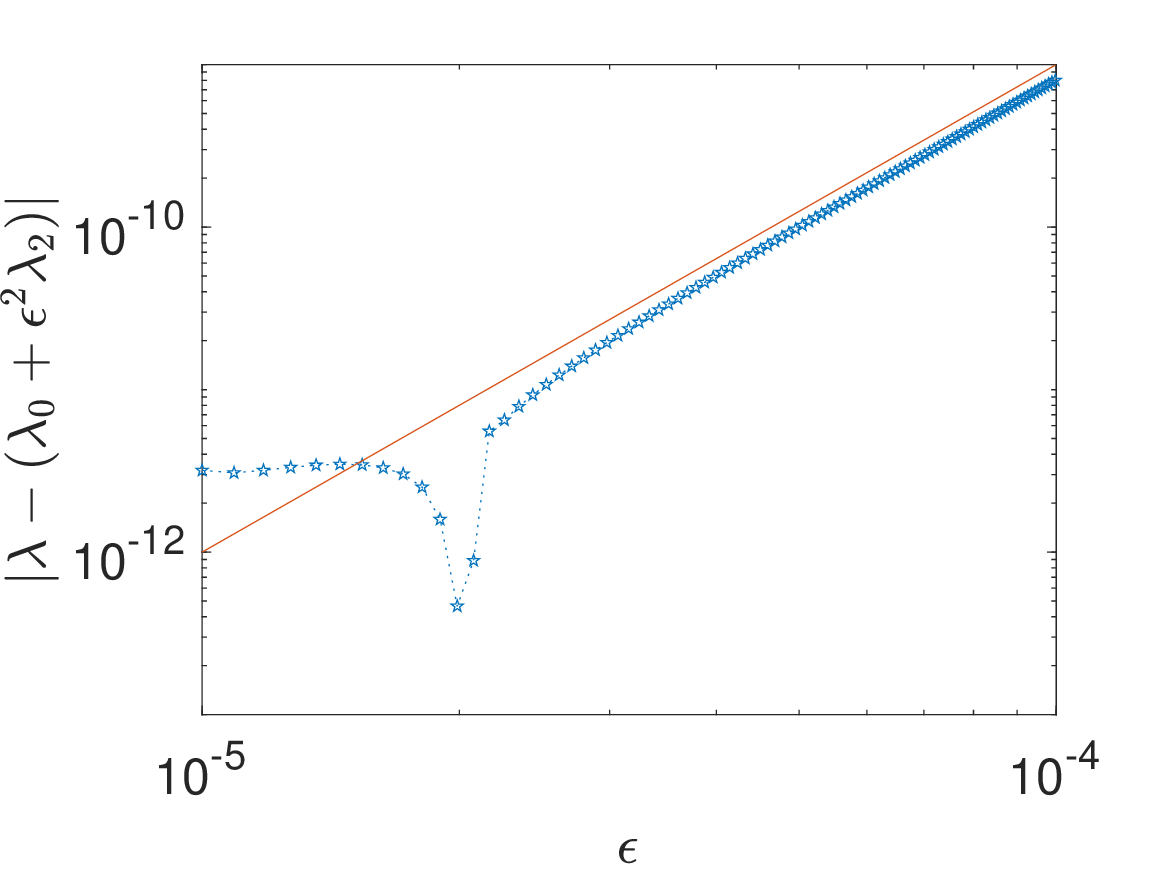}}
\caption{\it \textbf{Left:} An isola of a quartet instability in the deep-water gravity-capillary Whitham equation near $(\lambda_0,p_0)=(0.2177i,0.1363)$  with $\sigma=0.25$ at $\varepsilon=10^{-4}$. The circles are computed with a quasi-Newton iteration; the curves are the second order asymptotic prediction of \eqref{QuartetLam2}. \textbf{Center:} The imaginary part of the most unstable eigenvalue of \eqref{QuartetLam2} (solid line) against numerical computations from a quasi-Newton iteration (circles).   \textbf{Right:}  The error between the numerically computed and asymptotics \eqref{QuartetLam2} (starred line); the solid line marks $y=10^3\varepsilon^3$. \label{WhithamQuartet}}
\end{figure}

Consider a quartet resonance between $k_2$ and $k_1$ for $k_2 > k_1$ without loss of generality. Then, $k_2 = k_1 + 2$, and the solvability conditions \eqref{MatrixSolvabilityQuartet} simplify to 
\begin{equation}\label{MatrixSolvabilityResQuartet}
 \left(\begin{array}{cc} -1 & -i\tilde{c}_g(k_1+p_0)\\ -\beta & -i\tilde{c}_g(k_2+p_0)\beta\end{array}\right)\left(\begin{array}{c}\lambda_2\\ p_2\end{array}\right)=\left(\begin{array}{c} i(k_1+p_0)\left\{\alpha_{1}^{+}+\alpha_{1}^{-}-c_2+(\alpha_{2}^{-}+D)\beta\right\}\\  i(k_2+p_0)\left\{(\alpha_{2}^{+}+\alpha_{2}^{-}-c_2)\beta+(\alpha_{1}^{+}+D)\right\} \end{array}\right),
\end{equation} 
with solutions
\[ \lambda_2=i\frac{A\beta^2+B\beta+C}{\beta}\quad \textrm{and} \quad p_2=\frac{E\beta^2+F\beta+G}{\beta}.\]
Here, $A,B,C,E,F,$ and $G$ are real-valued coefficients given by
\begin{eqnarray*}
A&=&\frac{-\tilde{c}_g(k_2+p_0)(k_1+p_0)(\alpha_{2}^{-}+D)}{\tilde{c}_g(k_2+p_0)-\tilde{c}_g(k_1)},\\
B&=&\frac{-\tilde{c}_g(k_2+p_0)(k_1+p_0)(\alpha_{1}^{+}+\alpha_{1}^{-}-c_2)+\tilde{c}_g(k_1+p_0) (k_2+p_0)(\alpha_{2}^{+}+\alpha_{2}^{-}-c_2)}{\tilde{c}_g(k_2+p_0)-\tilde{c}_g(k_1+p_0)},\\
C&=&-\frac{-\tilde{c}_g(k_1+p_0) (k_2+p_0)(\alpha_{1}^{+}+D)}{\tilde{c}_g(k_2+p_0)-\tilde{c}_g(k_1+p_0)},\\
E&=&\frac{(k_1+p_0)(\alpha_{2}^{-}+D)}{\tilde{c}_g(k_2+p_0)-\tilde{c}_g(k_1+p_0)},\\
F&=&\frac{(k_1+p_0)(\alpha_{1}^{+}+\alpha_{1}^{-}-c_2)-(k_2+p_0)(\alpha_2^{+}+\alpha_{2}^{-}-c_2)}{\tilde{c}_g(k_2+p_0)-\tilde{c}_g(k_1+p_0)},\\
G&=&-\frac{(k_2+p_0)(\alpha_1^{+}+D)}{\tilde{c}_g(k_2+p_0)-\tilde{c}_g(k_1+p_0)}.
\end{eqnarray*} 

 If $\beta \in \mathbb{R}$, then $\lambda_2$ is imaginary, similar to what occurs for $\lambda_1$ in the triad case. If $\beta$ is complex, then $\beta=\rho e^{i\theta}$ and
\[ p_2=E\rho e^{i\theta}+F+G\rho^{-1}e^{-i\theta}. \]
It follows that 
\[ \text{Im}(p_2)=E\rho \sin(\theta)-G\rho^{-1}\sin(\theta). \]
To ensure the Floquet exponent remains real-valued, we must enforce either $\sin(\theta)=0$ (leading to $\beta \in \mathbb{R}$) or, more interestingly, $\rho^2=\frac{G}{E}$ and $\theta$ is free. For quartet high-frequency instabilities to occur, we must have $G/E > 0$. As in the triad case, this sign condition is exactly the Krein signature of the repeated eigenvalue $\lambda_0$. Our work shows that this Krein condition is necessary and sufficient to develop quartet high-frequency instabilities. 

Given the choice of $\beta$ above and assuming the Krein condition holds, the formula for the eigenvalue corrections is
\begin{equation}\label{QuartetLam2}
 \lambda_2=i\left( A\sqrt{\frac{G}{E}}(\cos(\theta)+i\sin(\theta)) +B+C\sqrt{\frac{E}{G}}(\cos(\theta)-i\sin(\theta)) \right), 
 \end{equation}
which has real and imaginary parts given by
\begin{subequations}
\begin{align}
\text{Re}(\lambda_2) &=\left(-A\sqrt{\frac{G}{E}}+C\sqrt{\frac{E}{G}}\right)\sin(\theta), \\
\text{Im}(\lambda_2)&=\left(A\sqrt{\frac{G}{E}}+C\sqrt{\frac{E}{G}}\right)\cos(\theta)+B.
\end{align}
\end{subequations}
respectively.
From these expressions, we see that quartet high-frequency instabilities also appear as ellipses in the complex spectral plane. However, unlike the triad case, these ellipses have semi-major axes that scale as $\mathcal{O}\left(\varepsilon^2\right)$, and their centers drift along the imaginary axis at a distance from $\lambda_0$ that grows like $\mathcal{O}\left(\varepsilon^2\right)$.

An example of a quartet high-frequency instability for the Kawahara equation appears in Figure \ref{KawaharaQuartet}, along with three different evaluations that demonstrate the validity of our asymptotic calculations. The left panel compares the shape and location of the asymptotic prediction of the isola (varying $\theta$ with fixed $\varepsilon$) to numerical computations of the spectrum computed with a quasi-Newton iteration.  The center panel compares plots of the imaginary component of the most unstable eigenvalue as a function of $\varepsilon$ according to \eqref{QuartetLam2} and to numerical computations. The right panel shows a log-log plot of the difference between the asymptotic and numerical predictions of the most unstable eigenvalue as a function of $\varepsilon$. From this log-log plot, we see that the observed errors are $\mathcal{O}(\varepsilon^3)$, as expected.  When $\varepsilon=10^{-3}$, errors are less than $10^{-8}$, hence the clear visible match between the asymptotic and numerical predictions in the left and center panel.

A quartet high-frequency instability isola for the deep water gravity-capillary Whitham equation appears in Figure \ref{WhithamQuartet}, along with similar evaluations supporting the validity of our calculations. Unlike the Kawahara equation, the numerics and asymptotics appear to disagree about the center of the isola. In the right panel of Figure \ref{WhithamQuartet}, a log-log plot of the difference between numerics and asymptotics shows that these errors have leading-order behavior $q\cdot\varepsilon^3$ for $q \approx 10^3$. That the exponent of this behavior is three gives confidence that our asymptotic calculations are correct. However, the prefactor $q$ is quite large.  As a consequence, we conclude that the cubic eigenvalue corrections of the Whitham equation are much larger than they were in the Kawahara example. This suggests that, if one were to compute the entire series, as in \cite{nicholls2007spectral}, the spectrum near this collision may have a small radius of convergence.

\subsubsection{Higher-Order Resonances}
If $|k_1-k_2|>2$, the solvability conditions simplify to
\begin{equation}\label{MatrixSolvabilityGenQuartet}
 \left(\begin{array}{cc} -1 & -i\tilde{c}_g(k_1+p_0)\\ -\beta & -i\tilde{c}_g(k_2+p_0)\beta\end{array}\right)\left(\begin{array}{c}\lambda_2\\ p_2\end{array}\right)=\left(\begin{array}{c} i(k_1+p_0)(\alpha_{1}^{+}+\alpha_{1}^{-}-c_2)\\  i(k_2+p_0)(\alpha_{2}^{+}+\alpha_{2}^{-}-c_2)\beta\end{array}\right).
\end{equation}
Note that $\beta$ appears in the second component of the right-hand side vector, as opposed to the first component when $|k_1 -k_2| = 2$. Assuming $\beta \neq 0$, we can divide $\beta$ from the second solvability condition, leaving behind solutions for $\lambda_2$ and $p_2$ that are independent of $\beta$:
\[ \lambda_2=\frac{i\left\{-(k_1+p_0)\tilde{c}_g(k_2+p_0))(\alpha_{1}^{+}+\alpha_{1}^{-}-c_2)+\tilde{c}_g(k_1+p_0))(k_2+p_0)(\alpha_{2}^{+}+\alpha_{2}^{-}-c_2)\right\}}{\tilde{c}_g(k_2+p_0)-\tilde{c}_g(k_1+p_0)}, \]
\[ p_2=\frac{(k_1+p_0)(\alpha_{1}^{+}+\alpha_{1}^{-}-c_2)-(k_2+p_0)(\alpha_{2}^{-}+\alpha_{2}^{+}-c_2)}{ \tilde{c}_g(k_2+p_0)-\tilde{c}_g(k_1+p_0)} \]
Since $\alpha_{j,\pm}\in \mathbb{R}$, we must have $\lambda_2\in i\mathbb{R}$ and $p_2\in\mathbb{R}$ for all $\beta\ne 0$.  Thus, higher-order resonances do not lead to instability at $\mathcal{O}\left(\varepsilon^2\right).$ In fact, one would expect that the the $m$th resonant high-frequency instability has the first opportunity to generate instability at $\mathcal{O}\left(\varepsilon^m\right)$. This is precisely what occurs in the finite-depth water wave problem.

\section{The Benjamin-Feir Instability\label{BenjaminFeir}}

Unstable eigenvalues of the Benjamin-Feir instability bifurcate from the zero eigenvalue of the flat-state spectrum. A direct calculation shows this eigenvalue has corresponding Floquet exponent $p_0 = 0$ and the wavenumbers $k$ in its eigenspace satisfy the resonance condition 
\begin{align*}
\lambda_0 = ic_0k - i\omega(k) = 0,
\end{align*}
where $c_0 = \omega(1)$. If we assume 
\begin{enumerate}
\item[(i)] $\omega(k)$ is odd, 
\item[(ii)] $\displaystyle \inf_{k \in \mathbb{N}_0 \setminus \{1\}}|\omega(k)/k - \omega(1)| > 0$, and
\item[(iii)] $\omega(k)/k \in \mathcal{C}^\alpha(\mathbb{R})$ for $\alpha \geq 3$,
\end{enumerate}
as in \cite{maspero2024full}, then the only solutions of the resonance condition are $k_1 = -1$, $k_2 = 0$, and $k_3 = 1$. Hence, $\lambda_0 = 0$ has geometric multiplicity three with eigenspace 
\begin{align*}
v_0 = \beta_{0,-1}\exp(-ix) + \beta_{0,0} + \beta_{0,1}\exp(ix),
\end{align*}
for $\beta_{0,0}$ and $\beta_{0,\pm1} \in \mathbb{C}$.

To capture unstable eigenvalues when $0 < \varepsilon \ll 1$, we expand the spectral data ($\lambda$,$v$) of \eqref{eq:spectral-problem} as power series in $\varepsilon$, much like the high-frequency instabilities. However, to ease the burden of the calculations that follow, we will re-express the spectral problem in the form
\begin{align}
Tv = \lambda v,& \quad \textrm{where} \quad \label{specProb} \\
T :=
i\Big(c(p+D)v - \omega(p+D)v - &2(p+D)(uv)\Big) \quad \textrm{and} \quad D = -i\partial_x. \nonumber
\end{align}
We normalize the eigenfunction $v$ according to
\begin{align*}
\int_{-\pi}^{\pi}ve^{ix}dx = 1,
\end{align*}
so that $\beta_{0,-1} = 1$ and all subsequent corrections of $v$ exclude the $\exp(-ix)$ mode. 

In addition to the spectral data, we expand the Floquet exponent $p$ in terms of a renormalized series 
\begin{align*}
p = \varepsilon p_1\left(1+r(\varepsilon)\right), \quad \textrm{with} \quad r(\varepsilon) = \sum_{n=1}^{\infty}r_n\varepsilon^n.
\end{align*}
As will be seen, $p_1$ belongs to a real-valued interval that is symmetric about zero, while $r(\varepsilon)$ rescales the radius of this interval as $\varepsilon$ increases. For the triad and quartet instabilities of the previous sections, the interval of $p$ corresponding to the unstable eigenvalues can have a non-zero midpoint. Thus, a traditional power series expansion is appropriate to capture these instabilities. In contrast, the interval of $p$ parameterizing the unstable eigenvalues of the Benjamin-Feir instability is always symmetric about zero. The renormalized expansion above reflects this additional constraint. For more details, see \cite{creedon2023high}. 

In light of (\emph{iii}), we define the operators
\begin{align*}
T_n := \frac{1}{n!}\frac{d^n T}{d\varepsilon^n} \biggr|_{\varepsilon=0}, \quad \textrm{for} \quad n \in \{0,...,\alpha\}.
\end{align*}
It follows that the $\mathcal{O}(\varepsilon)$ terms of \eqref{specProb} take the form
\begin{align}
T_0v_1 = \left(\lambda_1-T_1\right)v_0, \label{ord1BF}
\end{align}
yielding an inhomogeneous linear equation for $v_1$. To ensure the solvability of \eqref{ord1BF}, the Fredholm alternative requires that the inhomogeneous terms are orthogonal to the nullspace of the adjoint of $T_0$. Since $T_0$ is skew self-adjoint, we can instead enforce orthogonality between the inhomogeneous terms and $v_0$, leading to the following solvability conditions:
\begin{subequations}
\begin{align}
\lambda_1+i\tilde{c}_g(-1)p_1  - i\beta_{0,0} &= 0, \label{solvCond1a} \\
\Big(\lambda_1 +i\tilde{c}_g(0)p_1  \Big)\beta_{0,0} &= 0, \label{solvCond1b}\\ 
\Big(\lambda_1+i\tilde{c}_g(1)p_1 \Big)\beta_{0,1} + i\beta_{0,0} &= 0. \label{solvCond1c}
\end{align}
\end{subequations}
To simplify the appearance of these conditions, we continue to use $\tilde{c}_g(k) := \omega'(k)-c_0$ as the linear group velocity of \eqref{mainModel} in a frame traveling at velocity $c_0$. Note by (\emph{i}) that $\tilde{c}_g(-k) = \tilde{c}_g(k)$ for all $k \in \mathbb{R}$.

System \eqref{solvCond1a}-\eqref{solvCond1c} admits many possible solutions, all of which imply that $\lambda_1$ is imaginary and, hence, no instability is observed at $\mathcal{O}\left(\varepsilon\right)$. Hindsight shows that there is a unique solution of these conditions that leads to eigenvalues with positive real part at $\mathcal{O}\left(\varepsilon^2\right)$. In particular, one should choose $\lambda_1 = -i\tilde{c}_g(1)p_1$ and $\beta_{0,0} = 0$, while $p_1$ and $\beta_{0,1}$ remain arbitrary at this order.

Having solved \eqref{solvCond1a}-\eqref{solvCond1c}, we return to \eqref{ord1BF}, invert $T_0$ against its range, and solve for $v_1$. Upon doing so, we find
\begin{align*}
v_1 = v_{1,-2}\exp(-2ix) + v_{1,2}\exp(2ix)  + \beta_{1,0} + \beta_{1,1}\exp(ix) ,
\end{align*}
where $\beta_{1,0}, \beta_{1,1} \in \mathbb{C}$, 
\begin{align}
v_{1,-2} =  \frac{1}{c_0-\frac{\omega(2)}{2}}, \quad \textrm{and} \quad v_{1,2} = \frac{\beta_{0,1}}{c_0-\frac{\omega(2)}{2}}.
\end{align}
Note that $v_{1,\pm 2}$ are well-defined by (\emph{ii}).

At $\mathcal{O}\left(\varepsilon^2
\right)$, \eqref{specProb} becomes
\begin{align*}
T_0v_2 = (\lambda_1 - T_1)v_1 + (\lambda_2 - T_2)v_0.
\end{align*}
Like the previous order, we require the inhomogeneous terms to satisfy the Fredholm alternative, leading us (after some work) to the following solvability conditions:
\begin{subequations}
\begin{align}
\beta_{1,0} &= \frac{1+\beta_{0,1}}{\omega'(1)-\omega'(0)}, \label{solvCond2a}\\
i\lambda_2 - \tilde{c}_g(1)r_1p_1  - \left(U+Vp_1^2 \right)&=U\beta_{0,1}, \label{solvCond2b} \\
\beta_{0,1}\Big(i\lambda_2 - \tilde{c}_g(1)r_1p_1 + U + Vp_1^2\Big) &= -U, \label{solvCond2c}
\end{align}
\end{subequations}
where 
\begin{align}
U = \frac{1}{\omega'(0)-\omega'(1)} + \frac{1}{\omega(2) - 2c_0} \quad \textrm{and} \quad  V = -\frac12\omega''(1).
\end{align}
By combining \eqref{solvCond2b} and \eqref{solvCond2c} to eliminate $p_1^2$ and then to eliminate $\lambda_2$, respectively, we find
\begin{align}
\lambda_2 &= -i\tilde{c}_g(1)r_1p_1 + i\frac{U}{2}\cdot \frac{(1-\beta_{0,1}^2)}{\beta_{0,1}}, \label{lambda2v2}\\
p_1^2 &= -\frac{U}{2V}\cdot\frac{(1+\beta_{0,1})^2}{\beta_{0,1}}.\label{p12v2}
\end{align}
 Since $p_1 \in \mathbb{R}$, equation \eqref{p12v2} requires
\begin{align}
\frac{(1+\beta_{0,1})^2}{\beta_{0,1}} \in \mathbb{R} \quad \implies \quad \frac{(1+\beta_{0,1})^2}{\beta_{0,1}} = \frac{(1+\overline{\beta}_{0,1})^2}{\overline{\beta}_{0,1}}, \label{betaEqn}
\end{align}
where the overbar indicates complex conjugation. Solving \eqref{betaEqn}, we find two possibilities: $\beta_{0,1} \in \mathbb{R} \setminus \{0\}$ or $|\beta_{0,1}| = 1$. We exclude the first possibility, as then $\lambda_2$ would be imaginary by \eqref{lambda2v2}. Hence, 
\begin{align*}
\beta_{0,1} = \exp(i\theta), \quad \textrm{for} \quad \theta \in [0,2\pi),
\end{align*}
and equations \eqref{lambda2v2} and \eqref{p12v2} become
\begin{align}
\lambda_2 &= \pm i\tilde{c}_g(1)r_1\sqrt{2\Delta_{\textrm{BF}}}\cos\left(\frac{\theta}{2} \right) + U\sin(\theta),\label{lambda2v3} \\
p_1^2 &= 2\Delta_{\textrm{BF}}\cos^2\left(\frac{\theta}{2} \right), \label{p12v3}
\end{align}
respectively, where
\begin{align}
\Delta_{\textrm{BF}} = -\frac{U}{V} = 2\left(\frac{\frac{1}{\omega'(0)-\omega'(1)} + \frac{1}{\omega(2) - 2c_0} }{\omega''(1)}\right).
\end{align}
In order for $p_1 \in \mathbb{R}$ for any $\theta$, we must have \begin{equation} \Delta_{\textrm{BF}} > 0, \end{equation} which is a sufficient condition for the development of the Benjamin-Feir instability. This condition is equivalent to those derived rigorously in \cite{maspero2024full,bronski2016modulational,hur2015modulationala,hur2015modulationalb}.

\begin{figure}[tp]
\centerline{\includegraphics[width=1.05\textwidth]{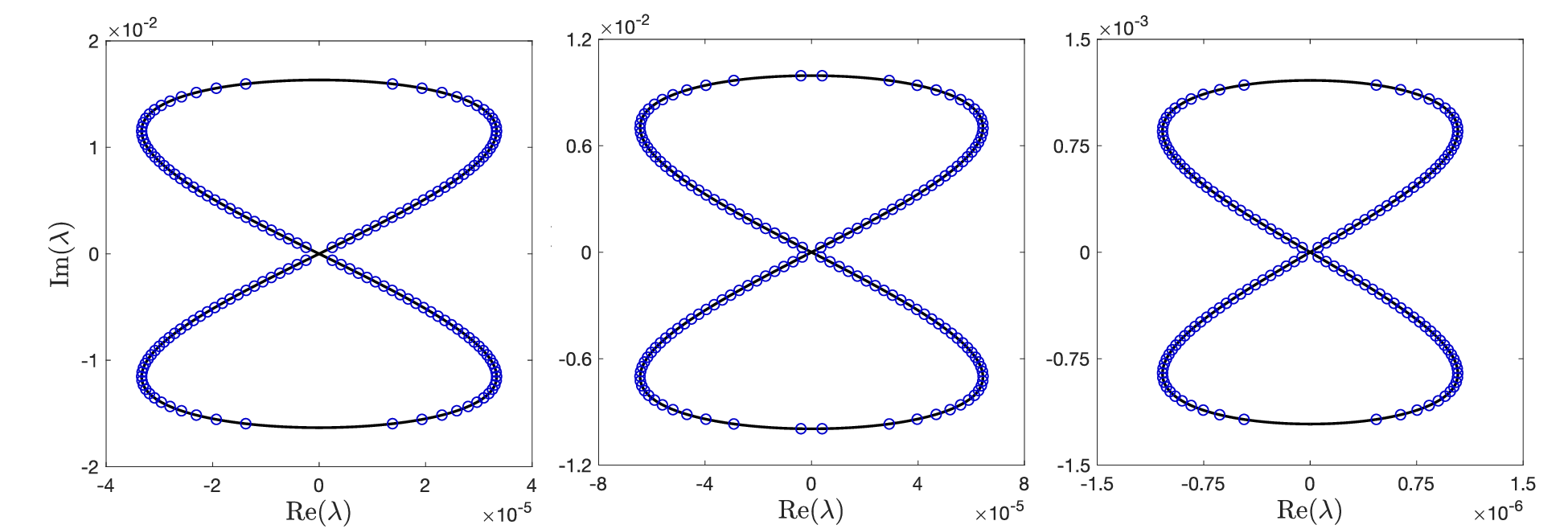}}
\caption{\it \textbf{Left:} Figure-eight curve for the Kawahara equation with $a = -3$, $b = 1$, and $\varepsilon = 10^{-2}$. The blue circles are computed using the Floquet-Fourier-Hill method; the black curves are the lemniscates predicted by \eqref{lemniscate}. \textbf{Center:} The same as the left plot but for the Whitham equation with $h = 2$. \textbf{Right:} The same as the left plot but for the capillary-Whitham equation with $h = 2$ and $\sigma = 3$.  \label{figureEightKawaharaWhitham}}
\end{figure}

Assuming $\Delta_{\textrm{BF}} > 0$, equation \eqref{p12v3} immediately implies
\begin{align}
-\sqrt{2\Delta_{\textrm{BF}}} \leq p_1 \leq \sqrt{2\Delta_{\textrm{BF}}}, \label{p1int}
\end{align}
so that $p_1$ is constrained to a real-valued interval that is symmetric about zero, as promised. The eigenvalues corresponding to these Floquet exponents lie on a closed curve of the form
\begin{subequations}
\begin{align}
\textrm{Im}\lambda &= \pm \varepsilon \tilde{c}_g(1)\sqrt{2\Delta_{\textrm{BF}}}\cos\left(\frac{\theta}{2}\right)\biggr(1 +  r_1 \varepsilon\biggr) + \mathcal{O}\left(\varepsilon^3\right), \\
\textrm{Re}\lambda &= \varepsilon^2U\Delta_{\textrm{BF}}\sin(\theta) + \mathcal{O}\left(\varepsilon^3\right).
\end{align} \label{lemniscate}
\end{subequations}
Truncating this parameterization after $\mathcal{O}\left(\varepsilon^2\right)$ and eliminating $\theta$ from the real and imaginary parts yields an equation for a lemniscate that approximates the figure-eight curve of Benjamin-Feir eigenvalues. Figure \ref{figureEightKawaharaWhitham} compares this lemniscate to eigenvalues computed numerically using the Floquet-Fourier-Hill method for the Kawahara, Whitham, and capillary-Whitham equations. Excellent agreement is found between the two computations of the eigenvalues, giving confidence in the correctness of our perturbation method.  \\

\noindent \textbf{Remark 4.1.} The parameter $r_1$ is still unknown at this order. Hence, there is remaining ambiguity in the imaginary component of the approximate lemniscate for the Benjamin-Feir eigenvalues. This ambiguity is resolved at the next order. In particular, at $\mathcal{O}\left(\varepsilon^3\right)$, the eigenvalue corrections are singular as $\theta \rightarrow 0$ and $\theta \rightarrow 2\pi$. The parameter $r_1$ is chosen to eliminate this singular behavior so that the third-order corrections to the Benjamin-Feir eigenvalues remain bounded for all $\theta$. This principle is known as the \emph{regular curve condition} and was first introduced in \cite{creedon2021higha,creedon2021highb,creedon2022high}. In our case, the regular curve condition yields $r_1 = 0$, similar to the water wave problem in finite depth \cite{creedon2023high}. In other instances, \emph{e.g.}, the water wave problem in infinite depth, $r_1 \neq 0$, which slightly alters the equation of the approximating lemniscate.  \\\\
\indent With \eqref{lemniscate} in hand, we can directly compare the relative strengths of the Benjamin-Feir and high-frequency instabilities \emph{across an entire family of weakly nonlinear dispersive models} for the first time. If the model admits triad high-frequency instabilities, these instabilities will always dominate the Benjamin-Feir instability at small amplitudes, as the triad instabilities grow as $\mathcal{O}\left(\varepsilon\right)$ while Benjamin-Feir grows as $\mathcal{O}\left(\varepsilon^2\right)$. In models for which there are no triad instabilities\footnote{This would occur when the Krein signature at the corresponding flat-state repeated eigenvalue is non-negative.}, then the Benjamin-Feir instability directly competes with the quartet high-frequency instability. From Section 3.0.3, the most unstable eigenvalue of the quartet instability has leading-order behavior $|-A\sqrt{G/E} + C\sqrt{E/G}|\varepsilon^2$, while  that of the Benjamin-Feir instability is $|U\Delta_{BF}|\varepsilon^2$. Thus, for Benjamin-Feir to overtake the quartet instability, one must have 
$$\left|-A\sqrt{G/E} + C\sqrt{E/G}\right| < \left|U\Delta_{BF}\right|.$$ 
Otherwise, Benjamin-Feir is subdominant to the quartet high-frequency instability. Should a model have neither triad nor quartet high-frequency instabilities, then Benjamin-Feir dominates at small amplitude by default. 

\subsection{The Akers--Milewski Equation}
From Table \ref{tabOfModels}, the dispersion relation of the Akers--Milewski equation with $\sigma=1$ is
\begin{align*}
\omega(k) = \textrm{sgn}(k)\big(1+|k|\big)^2,
\end{align*}
which violates assumption (\emph{iii}) in our investigation of the Benjamin-Feir instability. Nevertheless, we can modify our perturbation method to capture the Benjamin-Feir instability in this model as well. To the authors' knowledge, this is the first time that the Benjamin-Feir instability spectrum has been approximated analytically for a model with discontinuous dispersion.

Crucial to the modified perturbation method is a correct expansion of
\begin{align*}
\omega(p+k) = \omega\left(\varepsilon p_1(1+r(\varepsilon))+k\right),
\end{align*}
with respect to $\varepsilon$ for each $k \in \mathbb{Z}$. By inspection, one sees there are two distinct cases:
\begin{align*}
\omega(p+k) = \begin{cases} \textrm{sgn}(k)\big(1+|\varepsilon p_1(1+r(\varepsilon)) + k|\big)^2 & k \in \mathbb{Z}\setminus\{0\} \\ \textrm{sgn}(\varepsilon p_1)\big(1 + |\varepsilon p_1(1+r(\varepsilon))|\big)^2 & k = 0\end{cases},
\end{align*}
where we have used $|\varepsilon| \ll 1$ to simplify the arguments of the sign functions. Both cases can be simplified further if one restricts $\varepsilon > 0$ and uses the identity
\begin{align*}
|s+t| = |s| + \textrm{sgn}(s)t \quad \textrm{for} \quad |t| < |s|.
\end{align*}
Then, 
\begin{align}
\omega(p+k) = \begin{cases} \textrm{sgn}(k)\big(1+\varepsilon|p_1| + \textrm{sgn}(p_1)r(\varepsilon)\big)^2 & k \in \mathbb{Z}\setminus\{0\} \\ \textrm{sgn}(p_1)\big(1 + \varepsilon|p_1| + \textrm{sgn}(p_1)r(\varepsilon)\big)^2 & k = 0\end{cases}. \label{newOmegaExp}
\end{align}
In both cases, \eqref{newOmegaExp} is now a power series in $\varepsilon$, allowing us to proceed formally with the perturbation method.

\begin{figure}[tp]
\centerline{\includegraphics[width=0.7\textwidth]{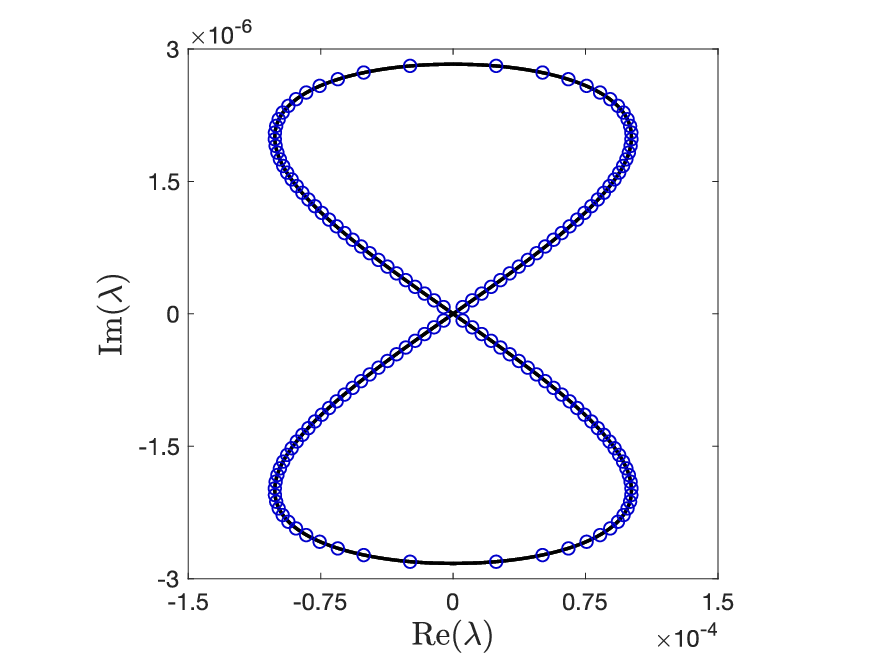}}
\caption{\it Figure-eight curve for the Akers--Milewski equation with $\varepsilon = 10^{-2}$ and $\sigma=1$. The blue circles are computed using the Floquet-Fourier-Hill method; the black curves are the lemniscates predicted by \eqref{lemniscate}.   \label{figureEightAkersMilewski}}
\end{figure}

As a consequence of \eqref{newOmegaExp}, the $k = 0$ Fourier mode is no longer in the nullspace\footnote{Thus, the geometric multiplicity of the zero eigenvalue drops from three to two for the Akers-Milewski equation, but the algebraic multiplicity remains equal to four.} of the operator $T_0$. Hence, at each order in $\varepsilon$ of the perturbation method, there are only two solvability conditions to satisfy, making the computations for the Akers--Milewski equation somewhat simpler. However, one also finds $\tilde{c}_g(1) = 0$ for this model, implying that the imaginary component of the figure-eight curve vanishes at first and second order in $\varepsilon$. Therefore, to get the leading-order lemniscate for the Akers--Milewski equation, one needs to carry out the perturbation method to $\mathcal{O}(\varepsilon^3)$. Working through the calculations, one indeed finds a lemniscate with $\mathcal{O}(\varepsilon^2)$ real component and $\mathcal{O}(\varepsilon^3)$ imaginary component. Figure \ref{figureEightAkersMilewski} compares this leading-order lemniscate with eigenvalues computed numerically by the Floquet-Fourier-Hill method, to excellent agreement.

\section{Conclusion\label{Conclusion}}

We have developed a unified perturbation framework for the spectral stability of
small-amplitude Stokes waves in weakly nonlinear model equations of the form
\eqref{mainModel}. By expanding the eigenvalues, eigenfunctions, and Floquet
exponent in the wave amplitude and treating the ratio $\beta$ of the colliding
flat-state eigenmodes as a free parameter, we organised the unstable spectrum
into sheets of spectral data in $(\textrm{Re}\lambda,\textrm{Im}\lambda,p)$–space.
For each fixed $\beta$ the solvability conditions determine the eigenvalue and
Floquet corrections, and varying $\beta$ traverses the full sheet associated
with a given collision. Slices of these sheets at fixed amplitude produce the isolas of instability discussed in \cite{akers2015modulational,creedon2021higha,creedon2021highb,creedon2022high}.

For non-zero Floquet collisions, we derived explicit asymptotics for triad and
quartet high-frequency instabilities. Triad resonances lead to elliptical
isolas of order $\mathcal{O}(\varepsilon)$ centred at the flat-state spectrum,
whereas quartet resonances produce $\mathcal{O}(\varepsilon^2)$ isolas whose
centres drift along the imaginary axis and whose semi-major axes also scale
like $\mathcal{O}(\varepsilon^2)$. Higher-order resonances do not generate
instability at these orders. At the triple collision at the origin, we
recovered the Benjamin–Feir (modulational) instability, obtained a lemniscate
approximation to the figure-eight spectrum, and formulated the perturbation
scheme so that it remains valid even when $\omega'(k)$ has jump
discontinuities. In particular, we derived, to our knowledge, the first analytic
approximation of the Benjamin–Feir spectrum for the Akers–Milewski equation,
showing that the method extends beyond smooth-dispersion models.

Because the same framework treats both high-frequency and Benjamin–Feir
instabilities across our family of models, it enables a direct comparison of
their growth rates in a common setting. When triad instabilities are present,
their $\mathcal{O}(\varepsilon)$ growth rates dominate the
$\mathcal{O}(\varepsilon^2)$ Benjamin–Feir growth; in the absence of triads, the
competition between quartet and Benjamin–Feir instabilities is governed by the
model-dependent coefficients appearing in their respective asymptotics. If a
model admits neither triad nor quartet high-frequency instabilities, the
Benjamin–Feir mechanism becomes the leading small-amplitude source of spectral
instability.

Throughout, asymptotic predictions were validated against
numerical spectra computed using the Floquet–Fourier–Hill method and a
quasi-Newton continuation scheme, with asymptotically correct agreement and
examples that highlight the influence of higher-order corrections. The
bifurcation structure revealed here is expected to inform future rigorous
existence results for unstable spectra and to guide further studies of nonlinear
dynamics near the onset of instability in dispersive wave models.

{\bf Disclaimer:} This report was prepared
as an account of work sponsored by an agency of the United States
Government. Neither the United States Government nor any agency
thereof, nor any of their employees, make any warranty, express or
implied, or assumes any legal liability or responsibility for the
accuracy, completeness, or usefulness of any information,
apparatus, product, or process disclosed, or represents that its
use would not infringe privately owned rights. Reference herein to
any specific commercial product, process, or service by trade name,
trademark, manufacturer, or otherwise does not necessarily
constitute or imply its endorsement, recommendation, or favoring
by the United States Government or any agency thereof. The views
and opinions of authors expressed herein do not necessarily state
or reflect those of the United States Government or any agency
thereof.\\

{\bf Declarations:} 
B. Akers was supported in part by the Air Force Office of Scientific Research (AFOSR) and the Joint Directed Energy Transition Office (DE-JTO) during the preparation of this manuscript.  R. Creedon acknowledges support from  the National Science Foundation, grant number NSF-DMS-2402044.

\bibliography{WiltonSheets.bib}{}

\end{document}